\overfullrule=0pt
\centerline {\bf On a minimax theorem: an improvement, a new proof and an overview of its applications}\par
\bigskip
\bigskip
\centerline {BIAGIO RICCERI}\par
\bigskip
\bigskip
{\bf Abstract.} Theorem 1 of [14], a minimax result for functions $f:X\times Y\to {\bf R}$, where $Y$ is a real interval,
was partially extended to the case where $Y$ is a convex set in a Hausdorff topological vector space ([15], Theorem 3.2).
In doing that, a key tool was a partial extension of the same result to the case where $Y$ is a convex set in ${\bf R}^n$
([7], Theorem 4.2). In the present paper, we first obtain a full extension of the result in [14] by means of a new proof fully
based on the use of the result itself via an inductive argument. Then, we present an overview of the various and numerous
applications of these results.\par
\bigskip
\bigskip
\bigskip
\bigskip
{\bf 1. Introduction}\par
\bigskip
In [14], we established the following result:\par
\medskip
THEOREM 1.A ([14], Theorem 1). - {\it Let $X$ be a topological space, $Y\subseteq {\bf R}$ an interval
and $f:X\times Y\to
{\bf R}$ a function satisfying the following conditions:\par
\noindent
$(a)$\hskip 5pt for each $y\in Y$, the function $f(\cdot,y)$ is lower semicontinuous and
inf-compact;\par
\noindent
$(b)$\hskip 5pt for each $x\in X$, the function $f(x,\cdot)$ is continuous and quasi-concave.\par
Then, at least one of the following assertions holds:\par
\noindent
$(i)$\hskip 5pt there exists $\hat y\in Y$ such that the function $f(\cdot,\hat y)$ has at least two global
minima;\par
\noindent
$(ii)$\hskip 5pt one has
$$\sup_Y\inf_Xf=\inf_X\sup_Yf\ .$$}\par
\medskip
Actually, in [14], $Y$ is assumed to be open. However, the same identical proof works for any interval $Y$ (see Remark 2.1 below).
\smallskip
Later, in [7], S.J.N. Mosconi obtained\par
\medskip
THEOREM 1.B ([7], Theorem 4.2). - {\it Let $X$ be a topological space, $Y\subseteq {\bf R}^n$ a non-empty convex set
and $f:X\times Y\to
{\bf R}$ a function satisfying the following conditions:\par
\noindent
$(a)$\hskip 5pt for each $y\in Y$, the function $f(\cdot,y)$ is lower semicontinuous and
inf-compact;\par
\noindent
$(b)$\hskip 5pt for each $x\in X$, the function $f(x,\cdot)$ is upper semicontinous and concave.\par
Then, the conclusion of Theorem 1.A holds.}\par
\medskip
Finally, in [15], using Theorem 1.B, we obtained\par
\medskip
THEOREM 1.C ([15], Theorem 3.2). - 
{\it Let $X$ be a topological space, $E$ a Hausdorff topological vector
space, $Y\subseteq E$ a non-empty convex set 
and $f:X\times Y\to
{\bf R}$ a function satisfying the following conditions:\par
\noindent
$(a)$\hskip 5pt for each $y\in Y$, the function $f(\cdot,y)$ is lower semicontinuous and
inf-compact;\par
\noindent
$(b)$\hskip 5pt for each $x\in X$, the function $f(x,\cdot)$ is upper semicontinous and concave.\par
Then, the conclusion of Theorem 1.A holds.}\par
\medskip
In comparing the above results, two natural questions arise: does Theorem 1.A hold if ``continuous" is relaxed
to ``upper semicontinuous" ? ; does Theorem 1.A hold if $Y$ is any non-empty convex set 
in a Hausdorff topological vector space ?\par
\medskip
The answer to the first question is negative. In this connection, consider the following\par
\medskip
EXAMPLE 1.1. - Let $X=\{x_0,x_1\}$ (with $x_0\neq x_1$ and $X$ equipped with the discrete topology) and let 
$f:X\times [0,1]\to {\bf R}$ be defined by
$$f(x_i,y)=\cases {f(x_i,y)=y & if $i=0, y\in [0,1]$\cr & \cr
f(x_i,y)=-y & if $i=1, y\in ]0,1]$\cr & \cr
f(x_i,0)=1 & if $i=1$\ .\cr}$$
Of course, $x_1$ is the only global minimum of $f(\cdot,y)$ for all $y\in ]0,1]$, while $x_0$ is the only global minimum of
$f(\cdot,0)$. Moreover, $f(x_i,\cdot)$ is upper semicontinuous and quasi-concave for $i=0,1$. However, we have
$$\sup_{[0,1]}\inf_Xf=0<1=\inf_X\sup_{[0,1]}f\ .$$\par
\medskip
To the contrary, the answer to the second question is positive. Indeed, we will prove\par
\medskip
THEOREM 1.1. - {\it Let $X$ be a topological space, $E$ a topological vector
space, $Y\subseteq E$ a non-empty convex set 
and $f:X\times Y\to
{\bf R}$ a function satisfying the following conditions:\par
\noindent
$(a)$\hskip 5pt for each $y\in Y$, the function $f(\cdot,y)$ is lower semicontinuous and
inf-compact;\par
\noindent
$(b)$\hskip 5pt for each $x\in X$, the function $f(x,\cdot)$ is continuous and quasi-concave.\par
Then, the conclusion of Theorem 1.A holds.}\par
\medskip
The aim of the present paper is twofold.\par
\smallskip
 On the one hand, we just wish to prove Theorem 1.1. We stress
that our proof of Theorem 1.1 is fully based on the use of Theorem 1.A, via an inductive argument.\par
\smallskip
In turn, using Theorem 1.1, we obtain\par
\medskip
THEOREM 1.2. - {\it Let $X$ be a topological space, $E$ a vector
space, $Y\subseteq E$ a non-empty convex set 
and $f:X\times Y\to
{\bf R}$ a function satisfying the following conditions:\par
\noindent
$(a)$\hskip 5pt for each $y\in Y$, the function $f(\cdot,y)$ is lower semicontinuous and
inf-compact;\par
\noindent
$(b)$\hskip 5pt for each $x\in X$, the function $f(x,\cdot)$ is concave.\par
Then, the conclusion of Theorem 1.A holds.}\par
\medskip
Hence, Theorem 1.2 is an improvement of Theorem 1.C obtained without resorting to Mosconi's result.\par
\smallskip
On the other hand, we wish to offer an overview of the several and various applications of Theorem 1.A (with its
``sequential" version) and Theorem 1.C known up to now ([12]-[21]).\par
\bigskip
{\bf 2. Proofs of Theorems 1.1 and 1.2}\par
\bigskip
As usual, a generic real-valued function $\varphi$ on a topological space $X$ is said to be inf-compact (resp.
inf-sequentially compact) if, for each $r\in {\bf R}$, the set $\varphi^{-1}(]-\infty,r])$ (called sub-level set) is compact (resp. 
sequentially compact). If $\varphi$ is defined on a convex set of a vector space, it is said to be quasi-concave if,
for each $r\in {\bf R}$, the set $\varphi^{-1}([r,+\infty[)$ is convex.\par
\smallskip
For each $n\in {\bf N}$, we put
$$S_n=\{(\lambda_1,...,\lambda_n)\in ([0,+\infty[)^n : \lambda_1+...+\lambda_n=1\}\ .$$
The core of our proof of Theorem 1.1 is to prove it first in the case where $Y=S_n$:\par
\medskip
LEMMA 2.1. - {\it Let $X$ be a topological space and let
$f:X\times S_n\to {\bf R}$ be a function satisfying the following conditions:\par
\noindent
$(a)$\hskip 5pt for each $y\in S_n$, the function $f(\cdot,y)$ is lower semicontinuous, inf-compact and has a unique
global minimum\ ;\par
\noindent
$(b)$\hskip 5pt for each $x\in X$, the function $f(x,\cdot)$ is continuous and quasi-concave.\par
Then, one has
$$\sup_{S_n}\inf_Xf=\inf_X\sup_{S_n}f\ .$$}\par
\smallskip
PROOF. We prove the theorem by induction on $n$. Clearly, it (trivially) holds for $n=1$. Now, assume that it is true
for $n=k$ ($k\geq 2$). We are going to prove that it is true for $n=k+1$. So, we are assuming that $f:X\times S_{k+1}\to
{\bf R}$ is a function satisfying $(a)$ and $(b)$ with $n=k+1$. 
Let $\psi:S_k\times [0,1]\to S_{k+1}$ be the continuous
function defined
by
$$\psi(\lambda_1,...,\lambda_{k},\mu)=(\mu\lambda_1,...,\mu\lambda_k,1-\mu)$$
for all $(\lambda_1,...,\lambda_k,\mu)\in S_k\times [0,1]$. 
Now, consider the function $\tilde f:X\times S_k\times [0,1]\to {\bf R}$ defined by
$$\tilde f(x,\lambda_1,...,\lambda_k,\mu)=f(x,\psi(\lambda_1,...,\lambda_k,\mu))$$
for all $(x,\lambda_1,...,\lambda_k,\mu)\in X\times S_k\times [0,1]$.
For each  $\mu\in [0,1]$ and for each $x\in X$, since $f(x,\cdot)$ is quasi-concave in $S_{k+1}$ and $\psi(\cdot,\mu)$ is
affine in $S_k$, it clearly follows that $\tilde f(x,\cdot,\mu)$ is quasi-concave in $S_k$. Therefore,
 by the induction assumption, we have
$$\sup_{y\in S_k}\inf_{x\in X}\tilde f(x,y,\mu) =\inf_{x\in X}\sup_{y\in S_k}\tilde f(x,y,\mu)\ .\eqno{(2.1)}$$
From $(2.1)$, we then infer
$$\sup_{(y,\mu)\in S_k\times [0,1]}\inf_{x\in X}\tilde f(x,y,\mu)=\sup_{\mu\in [0,1]}\sup_{y\in S_k}\inf_{x\in X}\tilde f(x,y,\mu)=
\sup_{\mu\in [0,1]}\inf_{x\in X}\sup_{y\in S_k}\tilde f(x,y,\mu)\ .\eqno{(2.2)}$$
Now, consider the function $g:X\times [0,1]\to {\bf R}$ defined by
$$g(x,\mu)=\sup_{y\in S_k}\tilde f(x,y,\mu)$$
for all $(x,\mu)\in X\times [0,1]$. Fix $\mu\in [0,1]$.
 From $(2.1)$, by compactness and semicontinuity, we infer the existence of
a point $(\hat x,\hat y)\in X\times S_k$ such that
$$\sup_{y\in S_k}\tilde f(\hat x,y,\mu)=\tilde f(\hat x,\hat y,\mu)=\inf_{x\in X}\tilde f(x,\hat y,\mu)\ .\eqno{(2.3)}$$
Now, let $x\in X$, with $x\neq \hat x$. By $(a)$ and $(2.3)$, we have
$$g(\hat x,\mu)=\tilde f(\hat x,\hat y,\mu)<\tilde f(x,\hat y,\mu)\leq g(x,\mu)\ .$$
In other words, $\hat x$ is the only global minimum of the function $g(\cdot,\mu)$ which is also lower semicontinuous
and inf-compact. Now, fix $x\in X$ and $r\in {\bf R}$. Set
$$C=\{u\in S_{k+1} : f(x,u)\geq r\}\ .$$
Of course, we have
$$\{\mu\in [0,1]: g(x,\mu)\geq r\}=\bigcup_{y\in S_k}\{\mu\in [0,1] : \tilde f(x,y,\mu)\geq r\}\ .\eqno{(2.4)}$$
Note that the right-hand side of $(2.4)$ is equal to the projection of the set
$\psi^{-1}(C)$ on $[0,1]$. But, for each $(\lambda_1,...,\lambda_{k+1})\in S_{k+1}$, 
we have
$$\psi^{-1}(\lambda_1,...,\lambda_{k+1})=\cases{\left \{\left ( {{\lambda_1}\over {1-\lambda_{k+1}}},...,{{\lambda_{k}}\over {1-\lambda_{k+1}}},
1-\lambda_{k+1}\right ) \right \} & if $\lambda_{k+1}\neq 1$\cr & \cr S_{k}\times \{0\} & if $\lambda_{k+1}=1$\ .\cr}$$
Hence, $\psi$ is onto $S_{k+1}$ and, by a classical result, for each compact and connected  set $D\subseteq S_{k+1}$, the set
$\psi^{-1}(D)$ is compact and connected. So, since $C$ is compact and connected (being convex), the set $\psi^{-1}(C)$ is connected and hence so is 
its projection on
$[0,1]$. Therefore, in view of $(2.4)$, the set $\{\mu\in [0,1]: g(x,\mu)\geq r\}$ is compact and connected. In other words, the
function $g(x,\cdot)$ is upper semicontinuous and quasi-concave in $[0,1]$. At this point, we can apply Theorem 1.A to $g$. So, we obtain
$$\sup_{\mu\in [0,1]}\inf_{x\in X}g(x,\mu)=\inf_{x\in X}\sup_{\mu\in [0,1]}g(x,\mu)\ .$$
Hence
$$\sup_{\mu\in [0,1}\inf_{x\in X}\sup_{y\in S_k}\tilde f(x,y,\mu)=\inf_{x\in X}\sup_{\mu\in [0,1]}\sup_{y\in S_k}\tilde f(x,y,\mu)=
\inf_{x\in X}\sup_{(y,\mu)\in S_k\times [0,1]}\tilde f(x,y,\mu)\ .\eqno{(2.5)}$$
Then, from $(2.2)$ and $(2.5)$, we get
$$\sup_{(y,\mu)\in S_k\times [0,1]}\inf_{x\in X}\tilde f(x,y,\mu)=\inf_{x\in X}\sup_{(y,\mu)\in S_k\times [0,1]}\tilde f(x,y,\mu)\ .
\eqno{(2.6)}$$
On the other hand, since $\psi(S_k\times [0,1])=S_{k+1}$, we have
$$\sup_{(y,\mu)\in S_k\times [0,1]}\inf_{x\in X}\tilde f(x,y,\mu)=\sup_{S_{k+1}}\inf_Xf$$
as well as
$$\inf_{x\in X}\sup_{(y,\mu)\in S_k\times [0,1]}\tilde f(x,y,\mu)=\inf_X\sup_{S_{k+1}}f$$
and so $(2.6)$ gives
$$\sup_{S_{k+1}}\inf_Xf=\inf_X\sup_{S_{k+1}}f\ ,$$
as claimed.\hfill $\bigtriangleup$\par
\medskip
A family of sets ${\cal C}$ is said to be filtering if for each pair $C_1, C_2\in {\cal C}$ there is $C_3\in {\cal C}$ such
that $C_1\cup C_2\subseteq C_3$.\par
\medskip
Now, we establish the following\par
\medskip
PROPOSITION 2.1. - {\it Let $X$ be a topological space, $Y$ a non-empty set, $y_0\in Y$ and $f:X\times Y\to {\bf R}$
a function such that $f(\cdot,y)$ is lower semicontinuous  for
all $y\in Y$ and inf-compact  for $y=y_0$. Assume also that there is a filtering
cover ${\cal C}$ of $Y$ such that
$$\sup_C\inf_Xf=\inf_X\sup_Cf$$
for all $C\in {\cal C}$.\par
Then, one has
$$\sup_Y\inf_Xf=\inf_X\sup_Yf\ .$$}\par
\smallskip
PROOF. Denote by ${\cal C}_0$ the family of all $C\in {\cal C}$ containing $y_0$. Clearly, ${\cal C}_0$ is a filtering
cover of $Y$. Arguing by contradiction, suppose that
$$\sup_Y\inf_Xf<\inf_X\sup_Yf\ .$$
Fix $r$ satisfying
$$\sup_Y\inf_Xf<r<\inf_X\sup_Yf$$
and, for each $C\in {\cal C}_0$, put
$$A_C=\left \{x\in X : \sup_{y\in C}f(x,y)\leq r\right \}\ .$$
Notice that $A_C\neq\emptyset$ since, otherwise, we would have 
$$r\leq \inf_X\sup_Cf=\sup_C\inf_Xf\leq \sup_Y\inf_Xf\ ,$$
against the choice of $r$. Now, observe that, if $C_1,...,C_k$ are finitely many members of ${\cal C}_0$, since ${\cal C}_0$ is
filtering, there is $\tilde C\in {\cal C}_0$ such that
$$\bigcup_{i=1}^kC_i\subseteq \tilde C\ .$$
This implies that
$$A_{\tilde C}\subseteq \bigcap_{i=1}^kA_{C_i}$$
and so $\bigcap_{i=1}^kA_{C_i}$ is non-empty. Therefore, $\{A_C\}_{C\in {\cal C}_0}$ is a family of closed subsets of the compact set 
$\{x\in X : f(x,y_0)\leq r\}$ possessing the finite intersection property. As as consequence, there would be $\tilde x\in
\cap_{C\in {\cal C}_0}A_C$. So,  since ${\cal C}_0$ is a cover of $Y$, we would have
$$\sup_Y\inf_Xf=\sup_{C\in {\cal C}_0}\sup_C\inf_Xf\leq\sup_{C\in {\cal C}_0}\sup_{y\in C}f(\tilde x,y)\leq r\ ,$$
against the choice of $r$.\hfill $\bigtriangleup$\par
The ``sequential" version of Proposition 2.1 is as follows:\par
\medskip
PROPOSITION 2.2. - {\it Let $X$ be a topological space, $Y$ a non-empty set, $y_0\in Y$ and $f:X\times Y\to {\bf R}$
a function such that $f(\cdot,y)$ is sequentially lower semicontinuous  for
all $y\in Y$ and inf-sequentially compact  for $y=y_0$. Assume also that there is an at most countable filtering
cover ${\cal C}$ of $Y$ such that
$$\sup_C\inf_Xf=\inf_X\sup_Cf$$
for all $C\in {\cal C}$.\par
Then, one has
$$\sup_Y\inf_Xf=\inf_X\sup_Yf\ .$$}\par
\smallskip
PROOF. Keep the notations of the above proof. An obvious inductive argument shows that there is a non-decreasing
sequence $\{C_k\}$ in ${\cal C}_0$ such that $\cup_{k\in {\bf N}}C_k=X$. So, $\{A_{C_k}\}$ turns out to be
a non-increasing sequence of non-empty sequentially closed subsets of the sequentially compact set $\{x\in X :
f(x,y_0)\leq r\}$. As a consequence, $\cap_{k\in {\bf N}}A_{C_k}\neq \emptyset$, and the proof goes as before.
\hfill $\bigtriangleup$\par
\medskip
{\it Proof of Theorem 1.1.} Denote by ${\cal P}$ the family of all convex polytopes of $Y$.  Of course, ${\cal P}$ is
a filtering cover of $Y$. Fix $P\in {\cal P}$. Let $x_1,...,x_n\in P$ be such that
$$P=\hbox {\rm conv}(\{x_1,...,x_n\})\ .$$
Consider the function $\eta : S_n\to P$ defined by
$$\eta(\lambda_1,...,\lambda_n)=\lambda_1x_1+...+\lambda_nx_n$$
for all $(\lambda_1,...,\lambda_n)\in S_n$. Plainly, the  function
$(x,\lambda_1,...,\lambda_n)\to f(x,\eta(\lambda_1,...,\lambda_n))$ satisfies in $X\times S_n$
the assumptions of Lemma 2.1, and so
$$\sup_{(\lambda_1,...,\lambda_n)\in S_n}\inf_{x\in X}f(x,\eta(\lambda_1,...,\lambda_n))=
\inf_{x\in X}\sup_{(\lambda_1,...,\lambda_n)\in S_n}f(x,\eta(\lambda_1,...,\lambda_n))\ .$$
Since $\eta(S_n)=P$, we then have
$$\sup_P\inf_Xf=\inf_X\sup_Pf\ .$$
Now, the conclusion follows from Proposition 2.1.\hfill $\bigtriangleup$\par
\medskip
{\it Proof of Theorem 1.2.}
Denote by ${\cal C}$ the family of all finite-dimensional convex subsets of $Y$. Fix $C\in {\cal C}$.
Denote by $L$ the linear span of $C$. Consider $L$ with the Euclidean topology. Since $C$ is convex, the relative interior
of $C$ (say $A$) is non-empty. By $(b)$, for each $x\in X$, the function $f(x,\cdot)_{|A}$ is continuous in $A$ and
one has
$$\sup_{y\in A}f(x,y)=\sup_{y\in C}f(x,y)\ .$$
By Theorem 1.1, we have
$$\sup_A\inf_Xf=\inf_X\sup_Af\ .$$
Therefore
$$\sup_A\inf_Xf\leq \sup_C\inf_Xf\leq\inf_X\sup_Cf=\inf_X\sup_Af$$
and so
$$\sup_C\inf_Xf=\inf_X\sup_Cf\ .$$
Now, the conclusion follows from Proposition 2.1\hfill $\bigtriangleup$\par
\medskip
REMARK 2.1. - As we said at the beginning, the proof of Theorem 1.A given in [14] holds for any interval
$Y$.  Actually, with the notation of [14], to get the lower semicontinuity of $\Psi$ in $X\times I$ it is enough to
apply Lemma 5 of [24] (which holds also when $f$ is quasi-concave in $I$). Furthermore, Theorem 1.A is still
true if, instead of $(a)$, we assume that, for each $y\in Y$, the function $f(\cdot,y)$ is sequentially lower semicontinuous
and inf-sequentially compact. In this case, to get the sequential lower semicontinuity of $\Psi$, it is enough to
apply Lemma 5 of [24] again, this time considering on $X$ the topology whose members are the sequentially open subsets of $X$.
\par
\bigskip
{\bf 3. A well-posedness theory}\par
\bigskip
In this section, we present a well-posedness theory for constrained minimization problems
which is based on the use of Theorem 1.A in its ``sequential" version (Remark 2.1).\par
\smallskip
In the sequel, $X$ is a  topological space, $J, \Phi$ are two
real-valued functions defined in $X$, and $a, b$ are two numbers in $[-\infty,+\infty]$,
with $a<b$. \par
\smallskip
If $a\in {\bf R}$ (resp. $b\in {\bf R}$), we denote by $M_a$ (resp. $M_b$) 
the set of all global minima of the function $J+a\Phi$ (resp. $J+b\Phi$), while
if $a=-\infty$ (resp. $b=+\infty$), $M_a$ (resp. $M_b$) stands for the empty set.
We adopt the conventions $\inf\emptyset=+\infty$, $\sup\emptyset=-\infty$.
\smallskip
We also set
$$\alpha:=\max\left \{ \inf_X \Phi,\sup_{M_b}\Phi\right \}\ ,$$
$$\beta:=\min\left \{ \sup_X \Phi,\inf_{M_a}\Phi\right \}\ .$$
Note that, by the next proposition, one has $\alpha\leq \beta$.\par
\smallskip
PROPOSITION 3.1. - {\it Let $Y$ be a nonempty set, 
$f, g:Y\to {\bf R}$ two functions, and $\lambda, \mu$ two
real numbers, with $\lambda<\mu$. Let $\hat y_{\lambda}$ 
 be a global minimum of the function $f+\lambda g$ and let
$\hat y_{\mu}$ be a global minimum of the function
$f+\mu g$.\par
Then, one has $$g(\hat y_{\mu})\leq g(\hat y_{\lambda})\ .$$ If either
$\hat y_{\lambda}$
or $\hat y_{\mu}$ is strict and $\hat y_{\lambda}\neq \hat y_{\mu}$, then
$$g(\hat y_{\mu})<g(\hat y_{\lambda})\ .$$}\par
\smallskip
PROOF. We have
$$f(\hat y_{\lambda})+\lambda g(\hat y_{\lambda})\leq f(\hat y_{\mu})+\lambda g(\hat y_{\mu})$$
as well as
$$f(y_{\hat \mu})+\mu g(\hat y_{\mu})\leq f(\hat y_{\lambda})+\mu g(\hat y_{\lambda})\ .$$
Summing, we get
$$\lambda g(\hat y_{\lambda})+\mu g(\hat y_{\mu})\leq +\lambda g(\hat y_{\mu})+\mu g(\hat y_{\lambda})$$
and so
$$(\lambda-\mu)g(\hat y_{\lambda})\leq (\lambda-\mu)g(\hat y_{\mu})$$
from which the first conclusion follows.
If either $\hat y_{\lambda}$ or $\hat y_{\mu}$ is strict and $\hat y_{\lambda}\neq \hat y_{\mu}$,
then one of the first two inequalities is strict and hence
so is the third one.
\hfill $\bigtriangleup$\par
\medskip
A usual, given a function $f:X\to {\bf R}$ and a set $C\subseteq X$, 
  we say that
the problem of minimizing $f$ over $C$ is well-posed if the following two
conditions hold:\par
\smallskip
\noindent
-\hskip 7pt the restriction of $f$ to $C$ has a unique global minimum,
say $\hat x$\ ;\par
\smallskip
\noindent
-\hskip 7pt every sequence $\{x_n\}$ in $C$ such that $\lim_{n\to \infty}f(x_n)=
\inf_Cf$, converges to $\hat x$.\par
\smallskip
 Clearly, when $f$ is inf-sequentially compact, the problem of minimizing
$f$ over a sequentially closed set
$C$ is well-posed if and only if $f_{|C}$ has a unique global
minimum.\par
\medskip
The basic result is as follows:\par
\medskip
THEOREM 3.1. - {\it 
Assume that $\alpha<\beta$ and that, for each $\lambda\in
]a,b[$, the function $J+\lambda\Phi$ is sequentially lower semicontinuous, inf-sequentially
compact and admits a unique global minimum in $X$.\par
Then, for each $r\in ]\alpha,\beta[$, 
the problem of minimizing $J$ over $\Phi^{-1}(r)$ is well-posed.\par
Moreover,
if we denote by $\hat x_r$  the unique global minimum
of $J_{|\Phi^{-1}(r)}$ $(r\in ]\alpha,\beta[)$,
the functions $r\to \hat x_r$ and $r\to J(\hat x_r)$ are continuous in
$]\alpha,\beta[$.}\par
\smallskip
PROOF. Fix $r\in ]\alpha,\beta[$
 and consider the function
$f:X\times {\bf R}\to {\bf R}$ defined by
$$f(x,\lambda)=J(x)+\lambda(\Phi(x)-r)$$
for all $(x,\lambda)\in X\times {\bf R}$. Clearly,
the the restriction of the function $f$ to $X\times
]a,b[$
satisfies all the assumptions of the variant of Theorem 1.A pointed out in Remark 2.1.
Consequently, since $(i)$ does not hold, we have
$$\sup_{\lambda\in ]a,b[}\inf_{x\in X}(J(x)+\lambda
(\Phi(x)-r))=
\inf_{x\in X}\sup_{\lambda\in ]a,b[} (J(x)+\lambda(\Phi(x)-r))
\ .\eqno{(3.1)}$$
Note that
$$\sup_{\lambda\in ]a,b[}\inf_{x\in X}f(x,\lambda)\leq
\sup_{\lambda\in [a,b]\cap {\bf R}}\inf_{x\in X}f(x,\lambda)\leq$$
$$\leq\inf_{x\in X}\sup_{\lambda\in [a,b]\cap {\bf R}}
f(x,\lambda)=
\inf_{x\in X}\sup_{\lambda\in ]a,b[}f(x,\lambda)$$
and so from $(3.1)$ it follows
$$\sup_{\lambda\in [a,b]\cap {\bf R}}\inf_{x\in X}(J(x)+\lambda
(\Phi(x)-r))=
\inf_{x\in X}\sup_{\lambda\in [a,b]\cap {\bf R}} (J(x)+\lambda(\Phi(x)-r))\ .\eqno{(3.2)}$$
Now, observe that the function $\inf_{x\in X}f(x,\cdot)$ is
upper semicontinuous in $[a,b]\cap {\bf R}$ and that
$$\lim_{\lambda\to +\infty}\inf_{x\in X}f(x,\lambda)=-\infty$$
if $b=+\infty$ (since $r>\inf_X \Phi$), and
$$\lim_{\lambda\to -\infty}\inf_{x\in X}f(x,\lambda)=-\infty$$
if $a=-\infty$ (since $r<\sup_X \Phi)$.
 From this, it clearly follows that
 there exists $\hat \lambda_r\in [a,b]\cap {\bf R}$ such that
$$\inf_{x\in X}f(x,\hat \lambda_r)=\sup_{\lambda\in [a,b]\cap {\bf R}}\inf_{x\in X}f(x,\hat \lambda_r)\ .$$
Since
$$\sup_{\lambda\in [a,b]\cap {\bf R}}f(x,\lambda)=\sup_{\lambda\in ]a,b[}f(x,\lambda)$$
for all $x\in X$, the sub-level sets of the function
$\sup_{\lambda\in [a,b]\cap {\bf R}}f(\cdot,\lambda)$ are sequentially compact.
 Hence, there exists $\hat x_r\in
X$ such that
$$\sup_{\lambda\in [a,b]\cap {\bf R}}f(\hat x_r,\lambda)=
\inf_{x\in X}\sup_{\lambda\in [a,b]\cap {\bf R}}f(x,\lambda)\ .$$
Then, thanks to $(3.2)$, $(\hat x_r,\hat \lambda_r)$ is a saddle-point of $f$, that is
$$J(\hat x_r)+\hat \lambda_r (\Phi(\hat x_r)-r)=\inf_{x\in X}(J(x)+\hat \lambda_r(\Phi(x)-r))=
J(\hat x_r)+\sup_{\lambda\in [a,b]\cap {\bf R}}\lambda(\Phi(\hat x_r)-r)\ .\eqno{(3.3)}$$
First of all, from $(3.3)$ it follows that $\hat x_r$ is a global minimum of
$J+\hat \lambda_r\Phi$.
We now show that $\Phi(\hat x_r)=r$. We
distinguish four cases.\par
\noindent
-\hskip 7pt $a=-\infty$ and $b=\infty$.  In this case, the equality $\Phi(\hat x_r)=
r$ follows from the fact that
$\sup_{\lambda\in {\bf R}}\lambda(\Phi(\hat x_r)-r)$ is finite.\par
\noindent
-\hskip 7pt $a>-\infty$ and $b=+\infty$. In this case, the finiteness of
$\sup_{\lambda\in [a,+\infty[}\lambda(\Phi(\hat x_r)-r)$ implies that $\Phi(\hat x_r)\leq
r$. But, if $\Phi(\hat x_r)<r$, from $(3.3)$, we would infer that
$\hat \lambda_r=a$ and so $\hat x_r\in M_a$. This would imply
 $\inf_{M_a}\Phi<r$, contrary to the choice of $r$.\par
\noindent
-\hskip 7pt $a=-\infty$ and $b<+\infty$. In this case, the finiteness of
$\sup_{\lambda\in ]-\infty,b]}\lambda(\Phi(\hat x_r)-r)$ implies that
$\Phi(\hat x_r)\geq r$. But, if $\Phi(\hat x_r)>r$, from $(3.3)$ again, we would infer
 $\hat \lambda_r=b$, and so $\hat x_r\in M_b$. Therefore, $\sup_{M_b}\Phi>r$, contrary to
the choice of $r$.\par
\noindent
-\hskip 7pt $-\infty<a$ and $b<+\infty$. In this case, if $\Phi(\hat x_r)\neq r$, as we
have just seen, we would have either $\inf_{M_a}\Phi<r$ or $\sup_{M_b}\Phi>r$, contrary
to the choice of $r$.\par
Having proved that $\Phi(\hat x_r)=r$, we also get that $\hat \lambda_r\in ]a,b[$.
Indeed, if $\hat \lambda_r\in \{a,b\}$, we would have either
$\hat x_r\in M_a$ or $\hat x_r\in M_b$ and so either $\inf_{M_a}\Phi
\leq r$ or $\sup_{M_b}\Phi\geq r$, contrary to the choice of $r$. From $(3.3)$ once again,
we furthermore infer that any global minimum of
$J_{|\Phi^{-1}(r)}$ (and $\hat x_r$ is so)
is a global minimum of $J+\hat \lambda_r\Phi$
in $X$. But, since $\hat \lambda_r\in ]a,b[$, $J+\hat \lambda_r\Phi$ has exactly one
global minimum in $X$ which, therefore, coincides with $\hat x_r$. Since
the sub-level sets of $J+\hat \lambda_r\Phi$ are sequentially compact, we then
conclude that any minimizing sequence in $X$ for $J+\hat \lambda_r\Phi$ converges
to $\hat x_r$. But any minimizing sequence in $\Phi^{-1}(r)$ for
$J$ is a minimizing sequence for
$J+\hat \lambda_r\Phi$, and so it converges to $\hat x_r$.
Consequently, the problem of minimizing $J$ over $\Phi^{-1}(r)$
is well-posed, as claimed. \par
Now, let us prove the other assertions made in thesis.
By Proposition 3.1, it clearly follows that the function $\lambda\to
\Phi(\hat y_{\lambda})$ is non-increasing in $]a,b[$ and that its
range is contained in $[\alpha,\beta]$. On the other hand, by the
first assertion of the thesis, this range contains $]\alpha,\beta[$.
Of course, from this it follows that the function $\lambda\to
\Phi(\hat y_{\lambda})$ is continuous in $]a,b[$. Now, observe
that the function $\lambda\to \inf_{x\in X}(J(x)+\lambda\Phi(x))$ is
concave and hence continuous in $]a,b[$. This, in particular, implies
that the function $\lambda\to J(\hat y_{\lambda})$ is continuous in
$]a,b[$.
Now, for each
$r\in ]\alpha,\beta[$, put
$$\Lambda_r=\{\lambda\in ]a,b[ : \Phi(\hat y_{\lambda})=r\}\ .$$
Let us prove that the multifunction $r\to \Lambda_r$ is upper semicontinuous in $]\alpha,\beta[$.
Of course, it is enough to show that the restriction of the multifunction to any bounded open
sub-interval of $]\alpha,\beta[$ is upper semicontinuous. So, let $s, t\in ]\alpha,\beta[$, with
$s<t$. Let $\mu, \nu\in ]a,b[$ be such that $\Phi(\hat y_{\mu})=t$, $\Phi(\hat y_{\nu})=s$. By
Proposition 3.1, we have 
$$\bigcup_{r\in ]s,t[}\Lambda_r\subseteq [\mu,\nu]\ .$$
Then, to show that the restriction of multifunction $r\to \Lambda_r$ to $]s,t[$
is upper semicontinuous, it is enough to prove that its graph is closed
in $]s,t[\times [\mu,\nu]$ ([6], Theorem 7.1.16).
But, this latter fact follows immediately from the continuity of the function
$\lambda\to \Phi(\hat y_{\lambda})$. At this point, we observe that, for
each $r\in ]\alpha,\beta[$, the function $\lambda\to \hat y_{\lambda}$ is
constant in $\Lambda_r$. Indeed, let $\lambda,\mu\in \Lambda_r$ with
$\lambda\neq \mu$. If it was $\hat y_{\lambda}\neq \hat y_{\mu}$, by
Proposition 3.1 it would follow
$$r=\Phi(\hat y_{\lambda})\neq \Phi(\hat y_{\mu})=r\ ,$$
an absurd. Hence, the function $r\to \hat x_r$, as composition of
the upper semicontinuous multifunction $r\to \Lambda_r$ and the continuous
function $\lambda\to \hat y_{\lambda}$, is continuous. Analogously,
the continuity of the function $r\to J(\hat x_r)$ follows observing that it
is the composition of $r\to \Lambda_r$ and the continuous function
$\lambda\to J(\hat y_{\lambda})$. The proof is complete.
\hfill $\bigtriangleup$\par
\medskip
REMARK 3.1. - It is important to remark that, under the assumptions of Theorem 3.1, we have
actually proved that, for each $r\in ]\alpha,\beta[$, there exists $\hat \lambda_r\in ]a,b[$
such that the unique global minimum of $J+\hat \lambda_r\Phi$ belongs to $\Phi^{-1}(r)$.\par
\medskip
When $a\geq 0$, we can obtain a conclusion dual to
that of Theorem 3.1,
 under the same key assumption.\par
\medskip
THEOREM 3.2 - {\it  Let $a\geq 0$.
 Assume that, for each $\lambda\in
]a,b[$, the function $J+\lambda\Phi$ is sequentially lower semicontinuous, inf-sequentially
compact and admits a unique global minimum in $X$.\par
Set
$$\gamma:=
\max\left \{ \inf_X J,\sup_{\hat M_a}J\right \}\ ,$$
$$\delta:=\min\left \{ \sup_X J,
\inf_{\hat M_b}J\right \}\ ,$$
where
$$\hat M_a=\cases {M_a & if $a>0$\cr & \cr \emptyset & if $a=0$\ ,\cr}$$
$$\hat M_b=\cases {M_b & if $b<+\infty$\cr & \cr \inf_X\Phi & if $b=+\infty$\ .\cr}$$
Assume that $\gamma<\delta$.\par
Then, for each $r\in ]\gamma,\delta[$, 
the problem of minimizing $\Phi$ over $J^{-1}(r)$ is well-posed.\par
Moreover,
if we denote by $\tilde x_r$  the unique global minimum
of $\Phi_{|J^{-1}(r)}$ $(r\in ]\gamma,\delta[)$,
the functions $r\to \tilde x_r$ and $r\to \Phi(\tilde x_r)$ are continuous in
$]\gamma,\delta[$.}\par
\smallskip
PROOF. Let $\mu\in ]b^{-1},a^{-1}[$. Then, since $\mu^{-1}\in ]a,b[$
and
$$\Phi+\mu J=\mu(J+\mu^{-1}\Phi)\ ,$$
we clearly have that the function $J+\mu\Phi$ has sequentially compact
sub-level sets and admits a unique global minimum. At this point, the conclusion
follows applying Theorem 3.1 with the roles of $J$ an $\Phi$ interchanged.
\hfill $\bigtriangleup$ \par
\medskip
We now state the version of Theorem 3.1 obtained in the setting of a
reflexive Banach space endowed with the weak topology.\par
\medskip
THEOREM 3.3. - {\it Let $X$ be a sequentially weakly closed set in a reflexive
real Banach space.
Assume that $\alpha<\beta$ and that, for each $\lambda\in
]a,b[$, the function $J+\lambda\Phi$ is sequentially weakly lower
semicontinuous, has bounded sub-level sets and has a unique global minimum in
$X$.  \par
Then, for each $r\in ]\alpha,\beta[$, 
the problem of minimizing $J$ over $\Phi^{-1}(r)$ is well-posed in the
weak topology. \par
Moreover,
if we denote by $\hat x_r$  the unique global minimum
of $J_{|\Phi^{-1}(r)}$ $(r\in ]\alpha,\beta[)$,
the functions $r\to \hat x_r$ and $r\to J(\hat x_r)$ are continuous in
$]\alpha,\beta[$, the first one in the weak topology.}
\smallskip
PROOF. Our assumptions clearly imply that,
for each $\lambda\in ]a,b[$, the sub-level sets of $J+\lambda\Phi$
are sequentially weakly compact, by the Eberlein-
Smulyan theorem. Hence, considering $X$ with the relative weak topology,
we are allowed to apply Theorem 3.1, from which the conclusion directly
follows.\hfill $\bigtriangleup$\par
\medskip
Analogously, from Theorem 3.2 we get\par
\medskip
THEOREM 3.4. - {\it Let $a\geq 0$ and let
 $X$ be a sequentially weakly closed set in a reflexive
real Banach space.
Assume that, for each $\lambda\in
]a,b[$, the function $J+\lambda\Phi$ is sequentially weakly lower
semicontinuous, has bounded sub-level sets and has a unique global minimum in
$X$. Assume also that $\gamma<\delta$, where $\gamma, \delta$ are defined as
in Theorem 3.2.\par
Then, for each $r\in ]\gamma,\delta[$, 
the problem of minimizing $\Phi$ over $J^{-1}(r)$ is well-posed in the
weak topology.
 \par
Moreover,
if we denote by $\tilde x_r$  the unique global minimum
of $\Phi_{|J^{-1}(r)}$ $(r\in ]\gamma,\delta[)$,
the functions $r\to \tilde x_r$ and $r\to \Phi(\tilde x_r)$ are continuous in
$]\gamma,\delta[$, the first one in the weak topology.}\par
\medskip
Finally, it is worth noticing that Theorem 3.1 also offers 
the perspective
of a novel way of seeing whether a given function possesses a global
minimum. Let us formalize this using Remark 3.1.\par
\medskip
THEOREM 3.5. - {\it
Assume that $b>0$ and that,
 for each $\lambda\in ]0,b[$,
the function $J+\lambda\Phi$ has sequentially compact sub-level sets and
 admits a unique global minimum, say
$\hat y_{\lambda}$. Assume also that
$$\lim_{\lambda\to 0^+}\Phi(\hat y_{\lambda})<\sup_X \Phi\ .\eqno{(3.4)}$$
Then, one has
$$ \lim_{\lambda\to 0^+}\Phi(\hat y_{\lambda})=\inf_M \Phi\ ,$$
where $M$ is the set of all global minima of $J$ in $X$.}\par
\smallskip
PROOF. We already know that the function $\lambda\to \Phi(\hat y_{\lambda})$
is non-increasing in $]a,b[$ and that its range is contained in
$[\alpha,\beta]$. We claim that
$$\beta=\lim_{\lambda\to 0^+}\Phi(\hat y_{\lambda})\ .$$
Assume the contrary.
 Let us apply Theorem 3.1, with
$a=0$ (so, $M_0=M$), using the conclusion pointed out in Remark 3.1. Choose $r$
satisfying
$$\lim_{\lambda\to 0^+}\Phi(\hat y_{\lambda})<r<\beta\ .$$
Then, (since also $\alpha<r$) it would exist $\hat \lambda_r\in ]0,b[$ such that
$\Phi(\hat y_{\hat\lambda_r})=r$,
 contrary to the choice of $r$.
 At this point, the conclusion follows directly from
 $(3.4)$.\hfill $\bigtriangleup$\par
\medskip
For the remainder of this section,  $X$ is an infinite-dimensional real Hilbert
space and $\Psi:X\to {\bf R}$ is a sequentially weakly continuous $C^1$
functional, with $\Psi(0)=0$.\par
\smallskip
For each $r>0$, set
$$S_r=\{x\in X : \|x\|^2=r\}$$
as well as
$$\gamma(r)=\sup_{x\in S_r}\Psi(x)\ .$$
Also, set
$$r^*=\inf\{r>0 : \gamma(r)>0\}\ .$$
In [25], M. Schechter and K. Tintarev developed a very elegant, transparent and precise theory
which can be summarized in the following result:\par
\medskip
THEOREM 3.A. - {\it Assume that $\Psi$
 has no local maximum in $X\setminus \{0\}$.
Moreover, let $I\subseteq ]r^*,+\infty[$ be an open interval such that,
for each $r\in I$, there exists a unique
$\hat x_r \in S_r$ satisfying $\Psi(\hat x_r)=\gamma(r)$. \par
Then, the following conclusions hold:\par
\noindent
$(i_1)$\hskip 5pt the function $r\to \hat x_r$ is continuous in $I$\ ;\par
\noindent
$(i_2)$\hskip 5pt the function $\gamma$ is $C^1$ and increasing in $I$\ ; \par
\noindent
$(i_3)$\hskip 5pt one has
$$\Psi'(\hat x_r)=2\gamma'(r)\hat x_r$$
for all $r\in I$.}
\medskip
The next result can
be regarded as the most complete fruit of a joint application of Theorems 3.A and 3.1.
\medskip
THEOREM 3.6. -  {\it 
Set
$$\rho=\limsup_{\|x\|\to +\infty}{{\Psi(x)}\over {\|x\|^2}}$$
and
$$\sigma=\sup_{x\in X\setminus \{0\}}{{\Psi(x)}\over {\|x\|^2}}\ .$$
Let $a, b$ satisfy
$$\max\{0,\rho\}\leq a<b\leq \sigma\ .$$
Assume that $\Psi$ has no local maximum in $X\setminus \{0\}$,
 and that, for each $\lambda\in ]a,b[$, the
functional 
$x\to \lambda\|x\|^2-\Psi(x)$
has a unique
global minimum, say $\hat y_{\lambda}$.
Let $M_a$ (resp. $M_b$ if $b<+\infty$ or
$M_b=\emptyset$ if $b=+\infty$) be
the set of all global minima of the functional
$x\to a\|x\|^2-\Psi(x)$ (resp. $x\to b\|x\|^2-\Psi(x)$
if $b<+\infty$).
Set
$$\alpha=\max\left \{0,\sup_{x\in M_b}\|x\|^2\right \} $$
and
$$\beta=\inf_{x\in M_a}\|x\|^2 \ .$$ 
Then, the following assertions hold:\par
\noindent
$(a_1)$\hskip 5pt one has $r^*\leq\alpha<\beta$\ ;\par
\noindent
$(a_2)$\hskip 5pt the function $\lambda\to g(\lambda):=
\|\hat y_{\lambda}\|^2$ is
decreasing in $]a,b[$ and its range is
$]\alpha,\beta[$\ ;\par
\noindent
$(a_3)$\hskip 5pt for each $r\in ]\alpha, \beta[$, the
point $\hat x_r:=\hat y_{g^{-1}(r)}$
is the unique global maximum of $\Psi_{|S_r}$ towards which
 every maximizing sequence in $S_r$ converges\ ;  \par
\noindent
$(a_4)$\hskip 5pt 
the function $r\to \hat x_r$ is continuous in $]\alpha,\beta[$\ ;\par
\noindent
$(a_5)$\hskip 5pt
the function $\gamma$ is $C^1$, increasing
 and strictly concave in
$]\alpha,\beta[$\ ; \par
\noindent
$(a_6)$\hskip 5pt
 one has
$$\Psi'(\hat x_r)=2\gamma'(r)\hat x_r$$
for all $r\in ]\alpha,\beta[$\ ;\par
\noindent
$(a_7)$\hskip 5pt one has
$$\gamma'(r)=g^{-1}(r)$$
for all $r\in ]\alpha,\beta[$.}\par
\smallskip
PROOF. First of all, observe that, by
Proposition 3.1, the function
  $g$ is non-increasing in $]a,b[$ and $g(]a,b[)\subseteq
[\alpha,\beta]$. 
Now, let $I\subset ]a,b[$ be a non-degenerate interval.
If $g$ was constant in $I$, then,
 by Proposition 3.1 again, 
 the function $\lambda\to \hat y_{\lambda}$
would be constant in $I$. Let $y^*$ be its unique value. Then,
$y^*$ would be a critical point of the functional $x\to
\lambda\|x\|^2-\Psi(x)$ for all $\lambda\in I$. That is to say
$$2\lambda y^*=\Psi'(y^*)$$
for all $\lambda\in I$. This would imply
that $y^*=0$,  and so (since $\Psi(0)=0$) we would have
$\inf_{x\in X}(\lambda\|x\|^2-\Psi(x))=0$ for all $\lambda\in I$, against the
fact that  $\inf_{x\in X}(\lambda\|x\|^2-\Psi(x))<0$ for all
$\lambda<\sigma$. Consequently, $g$ is decreasing in $]a,b[$, and so,
in particular, $\alpha<\beta$.  
Next, observe that
$$\lim_{\|x\|\to +\infty}(\lambda\|x\|^2-\Psi(x))=+\infty$$
for each $\lambda>\max\{0,\rho\}$. From this, 
recalling that $\Psi$ is sequentially weakly continuous, it 
clearly
follows that we can apply Theorem 3.1, taking
$J=-\Psi$ and $\Phi(\cdot)=\|\cdot\|^2$. 
Consequently (see Remark 3.1), for every $r\in ]\alpha,\beta[$, there
exists $\lambda_r\in ]a,b[$ such that $\|\hat y_{\lambda_r}\|^2=r$.
Therefore, by the strict monotonicity of $g$, we 
have $g(]a,b[)=]\alpha,\beta[$.
Now, let us prove $(a_3)$.
 Fix $r\in ]\alpha,\beta[$. Clearly, we have 
$$\|\hat x_r\|^2=r\ .$$
Since
$$g^{-1}(r)\|\hat x_r\|^2-\Psi(\hat x_r)\leq
g^{-1}(r)\|x\|^2-\Psi(x)$$
for all $x\in X$, we then have
$$\Psi(x)\leq \Psi(\hat x_r)$$
for all $x\in S_r$. Hence, $\hat x_r$ is a global maximum of
$\Psi_{|S_r}$. On the other hand, if $v$ is a global maximum 
of $\Psi_{|S_r}$,
 then 
$$g^{-1}(r)\|v\|^2-\Psi(v)=g^{-1}(r)\|\hat x_r\|^2-\Psi(\hat x_r)$$
and hence, since
$$\inf_{x\in X}(g^{-1}(r)\|x\|^2-\Psi(x))=g^{-1}(r)\|\hat x_r\|^2-
\Psi(\hat x_r)\ ,$$
we have $v=\hat x_r$. 
In other words, $\hat x_r$ is the unique global
maximum of $\Psi_{|S_r}$. Since the sub-level sets of the functional
$x\to g^{-1}(r)\|x\|^2-\Psi(x)$ are sequentially weakly compact, any minimizing sequence of this
functional in $X$ converges weakly to $\hat x_r$. Now, let
$\{w_n\}$ be any sequence in $S_r$ such that $\lim_{n\to \infty}
\Psi(w_n)=\gamma(r)$. Then, we have 
$$\lim_{n\to \infty}g^{-1}(r)\|w_n\|^2-\Psi(w_n)=
\inf_{x\in X}(g^{-1}(r)\|x\|^2-\Psi(x))$$
and so $\{w_n\}$ converges weakly to $\hat x_r$. But then, since $\lim_{n\to \infty}\|w_n\|=\|\hat x_r\|$
and $X$ is a Hilbert space, we have $\lim_{n\to \infty}\|w_n-x_r\|=0$ by a classical result.
Let us prove that $r^*\leq \alpha$. Arguing by contradiction,
assume that $\alpha<r^*$. Choose $r\in ]\alpha,\min\{r^*,\beta\}[$. 
Then, since $\gamma$ is non-decreasing in $]0,+\infty[$
(see Lemma 2.1 of [25]) and $\Psi$ is continuous,
we would have $\gamma(r)=0$, and so
$\Psi(\hat x_r)=0$, and this would contradict the fact that 
$\inf_{x\in X}(g^{-1}(r)\|x\|^2-\Psi(x))<0$ since $g^{-1}(r)<\sigma$.
At this point, we are allowed to apply Theorem 3.A taking
$I=]\alpha,\beta[$. Consequently, the function $\gamma$ is $C^1$ and
increasing in $]\alpha,\beta[$, and
  $(a_4)$, $(a_6)$ come directly from
$(i_1)$, $(i_3)$ respectively.
Fix $r\in ]\alpha,\beta[$ again. Since $\hat x_r$ is a critical point 
of the functional $x\to g^{-1}(r)\|x\|^2-\Psi(x)$, we have
$$2g^{-1}(r)\hat x_r=\Psi'(\hat x_r)$$
and then $(a_7)$ follows from a comparison
with $(a_6)$. Finally, from $(a_7)$, since
$g^{-1}$ is decreasing in $]\alpha,\beta[$, it follows that $\gamma$
is strictly concave there, and the proof is complete.\hfill
$\bigtriangleup$
\medskip
REMARK 3.2. - If the derivative of $\Psi$ is compact and if, for some
$\lambda>\rho$, the functional $x\to \lambda\|x\|^2-\Psi(x)$ 
 has at most two critical points in $X$, then the same functional
has a unique global minimum in $X$. Indeed, if this functional
had at least two global minima, taken into account that it satisfies the
classical Palais-Smale condition ([29], Example 38.25), it would have at
least three critical points by Corollary 1 of [9]. \par
\bigskip
{\bf 4. A strict minimax inequality theory}\par
\bigskip
In order to use the results of Section 1 to get the multiplicity of global minima, we need
to know that, in the considered case, the strict minimax inequality holds.\par
\smallskip
The present section is just devoted to a theory on this matter.\par
\smallskip
To state our results in a more compact form, we now fix some notations.\par
\medskip
Throughout this section, $X$ is a non-empty set, $\Lambda, Y$ are two topological
spaces, $y_0$ is a point in $Y$.\par
\smallskip
A family ${\cal N}$ of non-empty subsets of $X$ is said to be a weakly filtering
cover of $X$ if  for each $x_1, x_2\in X$
there is $A\in {\cal N}$ such that $x_1, x_2\in A$.
\smallskip
We denote by ${\cal G}$
the family of all lower semicontinuous functions $\varphi:Y\to [0,+\infty[$,
with $\varphi^{-1}(0)=\{y_0\}$, such that, for each neighbourhood $V$ of $y_0$,
one has
$$\inf_{Y\setminus V}\varphi>0\ .\eqno{(4.1)}$$
Moreover, we denote by ${\cal H}$ the family of all functions
$\Psi:X\times\Lambda\to Y$ such that, for each $x\in X$, 
$\Psi(x,\cdot)$ is continuous, injective, open, 
takes the value $y_0$ at a point $\lambda_x$ and the function
$x\to \lambda_x$ is not constant. 
Furthermore, we denote by ${\cal M}$ the family of all
functions $J:X\to {\bf R}$ whose set of all global minima
(noted by $M_{J}$) is non-empty.\par
\smallskip
Finally, for each $\varphi\in {\cal G}$, $\Psi\in {\cal H}$ 
 and $J\in {\cal M}$, we put
$$\theta(\varphi,\Psi,J)=\inf\left \{
{{J(x)-J(u)}\over {\varphi(\Psi(x,\lambda_u))}} : (u,x)\in M_{J}\times
X\hskip 3pt \hbox {\rm with}\hskip 3pt \lambda_x\neq \lambda_u\right
\}\ .$$
\bigskip
With such notations, our theory is summarized in the following result:\par
\medskip
THEOREM 4.1. - {\it Let $\varphi\in {\cal G}$, $\Psi\in {\cal H}$  
 and $J\in {\cal M}$. \par
Then, for each $\mu>\theta(\varphi,\Psi,J)$ and each weakly filtering cover
${\cal N}$ of $X$, there exists $A\in {\cal N}$ such that
$$\sup_{\lambda\in \Lambda}\inf_{x\in A}
(J(x)-\mu\varphi(\Psi(x,\lambda)))<
\inf_{x\in A}\sup_{z\in A}
(J(x)-\mu\varphi(\Psi(x,\lambda_z)))\ .$$}\par
\smallskip
PROOF. 
 Let $\mu>\theta(\varphi,\Psi,J)$ and let ${\cal N}$ be
a weakly filtering cover of $X$. 
 Choose $u\in M_{J}$ and $x_1\in X$,
with $\lambda_{x_1}\neq \lambda_u$, such that
$$J(x_1)-\mu\varphi(\Psi(x_1,\lambda_u))<J(u)\ .$$
Let $A\in
{\cal N}$ be such that $u, x_{1}\in A$. We have
$$0\leq\inf_{z\in A}\varphi(\Psi(x,\lambda_z))\leq
\varphi(\Psi(x,\lambda_x))=0$$
for all $x\in A$, and so, since $u$ is a global
minimum of $J$, it follows that
$$\inf_{x\in A}\sup_{z\in A}
(J(x)-\mu\varphi(\Psi(x,\lambda_z)))=\inf_{x\in A}
\left ( J(x)-\mu\inf_{z\in A}
\varphi(\Psi(x,\lambda_z))\right )$$
$$ =\inf_XJ=J(u)\ .\eqno{(4.2)}$$
Since the function
$\varphi(\Psi(x_1,\cdot))$ is lower semicontinuous at
$\lambda_u$, 
there are $\epsilon>0$ and
a neighbourhood $U$ of $\lambda_u$ such that
$$J(x_1)-\mu\varphi(\Psi(x_1,\lambda))<J(u)-\epsilon$$
for all $\lambda\in U$. 
So, we have
$$\sup_{\lambda\in U}\inf_{x\in A}
(J(x)-\mu\varphi(\Psi(x,\lambda)))\leq\sup_{\lambda\in U} 
 (J(x_1)-\mu\varphi(\Psi(x_1,\lambda)))\leq J(u)-\epsilon
\ .\eqno{(4.3)}$$
Since $\Psi(u,\cdot)$ is open,
 the set $\Psi(u,U)$ is a neighbourhood of $y_0$.
Hence, by $(4.1)$, we have
$$\nu:=\inf_{y\in Y\setminus \Psi(u,U)}\varphi(y)>0\ .\eqno{(4.4)}$$
Moreover, since $\Psi(u,\cdot)$ is injective,
 if $\lambda\not\in U$ then $\Psi(u,\lambda)\not
\in \Psi(u,U)$. So,
 from $(4.4)$, it follows that
$$\sup_{\lambda\in \Lambda\setminus U}\inf_{x\in A}
(J(x)-\mu\varphi(\Psi(x,\lambda))\leq J(u)-
\mu\inf_{\lambda\in \Lambda\setminus U}\varphi(\Psi(u,\lambda))\leq
J(u)-\mu\nu\ . \eqno{(4.5)}$$
Now, the conclusion comes directly from $(4.2)$, $(4.3)$, $(4.4)$ and
$(4.5)$.\hfill $\bigtriangleup$\par
\medskip
REMARK 4.1. - From the conclusion of Theorem 4.1 it clearly follows that,
for any set $D\subseteq\Lambda$ with $\lambda_x\in D$ for all $x\in A$,
one has
$$\sup_{\lambda\in D}\inf_{x\in A}
(J(x)-\mu\varphi(\Psi(x,\lambda)))<
\inf_{x\in A}\sup_{\lambda\in D}
(J(x)-\mu\varphi(\Psi(x,\lambda)))\ .$$
\medskip
REMARK 4.2. - From the definition of $\theta(\varphi,\Psi,J)$, it clearly follows that
$u\in M_J$ if and only if $u$ is a global minimum of the
function $x\to J(x)-\theta(\varphi,\Psi,J)\varphi(\Psi(x,\lambda_u))$.
So, when
$\theta(\varphi,\Psi,J)>0$, from knowing that
$$J(u)\leq J(x)$$
for all $x\in X$, we automatically get
$$J(u)\leq J(x)-\theta(\varphi,\Psi,J)\varphi(\Psi(x,\lambda_u))$$
for all $x\in X$, which is a much better inequality since
$\varphi(y)>0$ for all $y\in Y\setminus \{y_0\}$.\par 
\medskip
REMARK 4.3 - It is likewise important to observe that if
$\theta(\varphi,\Psi,J)>0$,
then the function $x\to \lambda_x$ is constant in $M_{J}$. As
a consequence, if $\theta(\varphi,\Psi,J)>0$ and
the function $x\to \lambda_x$ is injective,
then $J$ has a unique global minimum. In
particular, note that $x\to \lambda_x$ is injective
 when $\Psi(\cdot,\lambda)$ is injective for
all $\lambda\in \Lambda$.\par
\medskip
REMARK 4.4. - Remarks 4.2 and 4.3 show the interest in knowing
when $\theta(\varphi,\Psi,J)>0$. Theorem 4.1 can also be useful
for this. Indeed, if for some $\mu>0$, there is a weakly filtering cover
${\cal N}$ of $X$ such that 
$$\sup_{\lambda\in \Lambda}\inf_{x\in A}
(J(x)-\mu\varphi(\Psi(x,\lambda)))\geq
\inf_{x\in A}\sup_{z\in A}
(J(x)-\mu\varphi(\Psi(x,\lambda_z)))$$
for all $A\in {\cal N}$, then $\theta(\varphi,\Psi,J)\geq \mu$.\par
\medskip
Notice the following consequence of Theorem 4.1:\par
\medskip
THEOREM 4.2. - {\it Let $Y$ be a inner product space, and let
$I:X\to {\bf R}$, $\Phi:X\to Y$ and $\mu>0$ be such that
the function $x\to I(x)+\mu\|\Phi(x)\|^2$ has a global minimum.\par
Then, at least one of the following assertions holds:\par
\noindent
$(a)$\hskip 5pt for each weakly filtering cover ${\cal N}$ of $X$, there exists
$A\in {\cal N}$ such that
$$\sup_{\lambda\in Y}\inf_{x\in A}(I(x)+\mu(2\langle\Phi(x),\lambda\rangle-
\|\lambda\|^2))<
\inf_{x\in A}\sup_{\lambda\in \Phi(A)}(I(x)+\mu(2\langle\Phi(x),\lambda\rangle-\|\lambda\|^2))\ ;$$
\noindent
$(b)$\hskip 5pt for each global minimum $u$ of $x\to I(x)+\mu\|\Phi(x)\|^2$,
one has
$$I(u)\leq I(x)+2\mu(\langle\Phi(x),\Phi(u)\rangle-\|\Phi(u)\|^2)$$
for all $x\in X$.}\par
\smallskip
PROOF. Take $\Lambda=Y$, $y_0=0$. For each $x\in X$, $y, \lambda\in Y$, set
$$\varphi(y)=\|y\|^2\ ,$$
$$\Psi(x,\lambda)=\Phi(x)-\lambda$$
and
$$J(x)=I(x)+\mu\|\Phi(x)\|^2\ .$$
So that
$$J(x)-\mu\varphi(\Psi(x,\lambda))=
I(x)+\mu(2\langle\Phi(x),\lambda\rangle -\|\lambda\|^2)\ .$$
With these choices, 
 $(b)$ is equivalent to the inequality
$$\mu\leq \theta(\varphi,\Psi,J)\ .$$
Now, the conclusion is a direct consequence of Theorem 4.1. \hfill $\bigtriangleup$ \par
\medskip
In turn, from Theorem 4.2, we get\par
\medskip
THEOREM 4.3. -  {\it Let $X$ be a non-empty set, $x_0\in X$, $Y$ a real inner product space,
$I:X\to {\bf R}$, $\Phi:X\to Y$, with $I(x_0)=0$, $\Phi(x_0)=0$, and $\mu>0$. Assume that
$$\inf_{x\in X}I(x)<0\leq\inf_{x\in X}(I(x)+\mu\|\Phi(x)\|^2)\ .$$
Then,  for each weakly filtering cover ${\cal N}$ of $X$, there exists
$A\in {\cal N}$ such that
$$\sup_{y\in Y}\inf_{x\in A}(I(x)+\mu(2\langle\Phi(x),y\rangle-
\|y\|^2))<
\inf_{x\in A}\sup_{y\in \Phi(A)}(I(x)+\mu(2\langle\Phi(x),y\rangle-\|y\|^2))\ .$$}\par
\smallskip
PROOF. The assumptions imply that $x_0$ is a global minimum of $x\to I(x)+\mu\|\Phi(x)\|^2$. But, at the same time,
since $\inf_XI<0$, $x_0$ is not a global minimum of $I$. Hence, $(b)$ of Theorem 4.2 does not hold and so $(a)$ holds. \hfill
$\bigtriangleup$
\bigskip
{\bf 5. Multiplicity of global minima}\par
\bigskip
In this section, we apply the results stated in Section 1 to obtain multiple global minima. \par
\medskip
THEOREM 5.1. - {\it 
 Let $X$ be a topological space
 and $J, \Phi:X\to {\bf R}$ two functions satisfying
the following conditions:\par
\noindent
$(a_1)$\hskip 5pt for each $\lambda>0$, the
function $J+\lambda\Phi$ has compact and
closed sub-level sets\ ;\par
\noindent
$(b_1)$ there exist $\rho\in ]\inf_X\Phi,\sup_X\Phi[$ and
$u_1, u_2\in X$ such that
$$\Phi(u_1)<\rho<\Phi(u_2)$$
and
$${{J(u_1)-\inf_{\Phi^{-1}(]-\infty,\rho])}J}\over {\rho-\Phi(u_1)}}<
{{J(u_2)-\inf_{\Phi^{-1}(]-\infty,\rho])}J}\over {\rho-\Phi(u_2)}}\ .$$
Under such hypotheses, there exists $\lambda^*>0$ such that the function
$J+\lambda^*\Phi$ has at least two global minima.}\par
\smallskip
PROOF.  Observe that, in view of Theorem 1 of [1],
condition $(b_1)$ is equivalent to the
inequality
$$\sup_{\lambda\geq 0}\inf_{x\in X}(J(x)+\lambda(\Phi(x)-\rho))<
\inf_{x\in X}\sup_{\lambda\geq 0}(J(x)+\lambda(\Phi(x)-\rho))\ . $$
On the other hand, since the function $\lambda\to \inf_{x\in X}
(J(x)+\lambda(\Phi(x)-\rho))$ is concave (and real-valued) in $]0,+\infty[$,
it is lower semicontinuous in $[0,+\infty[$ and so
$$\sup_{\lambda\geq 0}\inf_{x\in X}(J(x)+\lambda(\Phi(x)-\rho))=
\sup_{\lambda> 0}\inf_{x\in X}(J(x)+\lambda(\Phi(x)-\rho))\ .$$
Consequently, condition $(b_1)$ is equivalent to the inequality
$$\sup_{\lambda>0}\inf_{x\in X}(J(x)+\lambda(\Phi(x)-\rho))<
\inf_{x\in X}\sup_{\lambda>0}(J(x)+\lambda(\Phi(x)-\rho))\ . $$
Now, we can apply Theorem 1.A taking $I=]0,+\infty[$ and
$$\Psi(x,\lambda)=J(x)+\lambda(\Phi(x)-\rho)\ ,$$
and the conclusion follows.\hfill $\bigtriangleup$
\medskip
A suitable application of Theorem 5.1 gives the following result:\par
\medskip
THEOREM 5.2. - {\it Let $S$ be a topological space and
$F,\Phi:S\to {\bf R}$ two lower semicontinuous
functions satisfying the following conditions:
\par
\noindent
$(a_2)$\hskip 5pt the function $\Phi$ is inf-compact\ ;\par
\noindent
$(b_2)$\hskip 5pt for some $a>0$, one has
$$\inf_{x\in \Phi^{-1}(]a,+\infty[)}{{F(x)}\over
{\Phi(x)}}=-\infty\ . $$
Under such hypotheses, for each $\rho$ large enough, there exists
$\lambda^*_{\rho}>0$ such that the restriction of the function
 $F+\lambda^*_{\rho}\Phi$ to $\Phi^{-1}(]-\infty,\rho])$
 has at least two global minima.} \par
\smallskip
PROOF. Fix
 $\rho_0>\inf_X\Phi$, $x_0\in \Phi^{-1}(]-\infty,\rho_0[)$ 
 and $\lambda$ satisfying
$$\lambda>
{{F(x_0)-\inf_{\Phi^{-1}(]-\infty,\rho_0])}F}\over
{\rho_0-\Phi(x_0)}}\ .$$
Hence, one has
$$F(x_0)+\lambda\Phi(x_0)<\lambda\rho_0+\inf_{\Phi^{-1}(]-\infty,\rho_0])}F
\ .\eqno{(5.1)}$$
Since $\Phi^{-1}(]-\infty,\rho_0])$ is compact, by lower
semicontinuity, there is
$\hat x\in \Phi^{-1}(]-\infty,\rho_0])$ such that
$$F(\hat x)+\lambda\Phi(\hat x)=
\inf_{x\in \Phi^{-1}(]-\infty,\rho_0])}(F(x)+\lambda\Phi(x))\ .\eqno{(5.2)}$$
We claim that $\Phi(\hat x)<\rho_0$. Arguing by contradiction,
assume that $\Phi(\hat x)\geq \rho_0$. Then, in view of $(5.1)$, we would have
$$F(x_0)+\lambda\Phi(x_0)<F(\hat x)+\lambda\Phi(\hat x)$$
against $(5.2)$. By $(b_2)$, there is a sequence $\{u_n\}$ in
$\Phi^{-1}(]a,+\infty[)$ such that 
$$\lim_{n\to \infty}{{F(u_n)}\over {\Phi(u_n)}}=-\infty\ .$$
Now, set
$$\gamma=\min\left \{ 0,\inf_{x\in \Phi^{-1}(]-\infty,\rho_0])}(F(x)+\lambda\Phi(x))\right \}$$
and fix $\hat n\in {\bf N}$ so that
$${{F(u_{\hat n})}\over {\Phi(u_{\hat n})}}<-\lambda + {{\gamma}\over {a}}\ .$$
We then have
$$F(u_{\hat n})+\lambda\Phi(u_{\hat n}) < {{\gamma}\over {a}}\Phi(u_{\hat n})\leq \gamma\ .$$
Hence, if we put
$$\rho^*=\Phi(u_{\hat n})\ ,$$
we have
$$\inf_{x\in \Phi^{-1}(]-\infty,\rho^*])}(F(x)+\lambda\Phi(x))<
\inf_{x\in \Phi^{-1}(]-\infty,\rho_0])}(F(x)+\lambda\Phi(x))\ .$$
At this point, for each $\rho\geq \rho^*$, we realize that it is possible
to apply Theorem 5.1 taking $X=\Phi^{-1}(]-\infty,\rho])$
and $J=F+\lambda\Phi$. Indeed, with these choices and
taking $u_1=\hat x$, $u_2=u_{\hat n}$,  the left-hand side of
the last inequality 
in $(b_1)$ is zero, while the right-hand side is positive. Consequently, there exists $\hat \lambda_{\rho}>0$ such that
 the restriction
of the function $F+\lambda\Phi+\hat\lambda_{\rho}\Phi$ to
$\Phi^{-1}(]-\infty,\rho])$ has at least two global minima. So,
 the conclusion follows taking $\lambda^*_{\rho}=\lambda+
\hat\lambda_{\rho}$.\hfill $\bigtriangleup$
\medskip
It is worth noticing the following consequence of Theorem 5.2.\par
\medskip
THEOREM 5.3. - {\it Let $S$ be a cone in a real vector space equipped with
a (not necessarily vector) topology and let $F, \Phi:S\to {\bf R}$
be two lower semicontinuous
functions satisfying the following conditions:
\par
\noindent
$(a_3)$\hskip 5pt the function $\Phi$ is
positively homogeneous of degree $\alpha$ and inf-compact ;\par
\noindent
$(b_3)$\hskip 5pt the function $F$ is positively homogeneous
of degree $\beta>\alpha$ and there is $\tilde x\in S$
such that $F(\tilde x)<0<\Phi(\tilde x)$\ .\par
Under such hypotheses, there exists $\rho^*>\inf_S\Phi$ such that
the restriction of the function $F+\Phi$ to $\Phi^{-1}(]-\infty,
\rho^*])$ has at least two global minima.}\par
\smallskip
PROOF. Clearly, we have
$$\lim_{\lambda\to +\infty}{{F(\lambda\tilde x)}\over {\Phi(\lambda\tilde x)}}=
\lim_{\lambda\to +\infty}{{F(\tilde x)}\over {\Phi(\tilde x)}}\lambda^{\beta-\alpha}=-\infty\ .$$
So, the hypotheses of Theorem 5.2 are satisfied and hence there exist
$\rho>\inf_S\Phi$ and $\lambda>0$ such that the
restriction of the function $F+\lambda\Phi$ to
$\Phi^{-1}(]-\infty,\rho])$ has at least two global
minima, say $v_1, v_2$. Now, observe that
$$\lambda^{\beta\over \alpha-\beta}(F(x)+\lambda\Phi(x))=
F(\lambda^{1\over \alpha-\beta} x)+\Phi(\lambda^{1\over \alpha-\beta} x)$$
for all $x\in S$. From this, it easily follows that the points
$\lambda^{1\over \alpha-\beta}v_1$ and $\lambda^{1\over \alpha-\beta}v_2$ are two global minima of the restriction of the function
$F+\Phi$ to $\Phi^{-1}(]-\infty,\lambda^{\alpha\over \alpha-\beta}\rho])$, that is the conclusion.\hfill $\bigtriangleup$
\medskip
REMARK 5.1. - We also remark that the number $\rho^*$ in the conclusion of Theorem 5.3 can be unique. In this connection, a very simple example
is provided by taking $S={\bf R}$, $\Phi(x)=x^2$ and $F(x)=-x^3$.
Actually, it is seen at once that, if $r>0$, the restriction of the function $x\to x^2-x^3$ to $[-r,r]$ has a unique global minimum when
$r\neq 1$ and exactly two global minima when $r=1$.\par
\medskip
With the notations of Section 4, 
a joint application of Theorem 1.2 and Theorem 4.1 gives\par
\medskip
THEOREM 5.4. - {\it Let $\varphi\in {\cal G}$, $\Psi\in {\cal H}$  
 and $J\in {\cal M}$. 
Moreover, assume that $X$ is a topological space, that
$\Lambda$ is a real vector space and that
$\varphi(\Psi(x,\cdot))$ is convex for each $x\in X$. Finally, let
$\mu>\theta(\varphi,\Psi,J)$ and ${\cal N}$ be a weakly filtering cover of
$X$
such that, for each $A\in {\cal N}$,
 the function $x\to J(x)-\mu\varphi(\Psi(x,\lambda))$ is
lower semicontinuous and inf-compact in $A$ for all $\lambda\in 
\hbox {\rm conv}(\{\lambda_x : x\in X\})$.
\par
Under such hypotheses,
there exist $A\in {\cal N}$ and $\lambda^*\in
\hbox {\rm conv}(\{\lambda_x : x\in A\})$
 such
that the restriction
of the function $x\to J(x)-\mu\varphi(\Psi(x,\lambda^*))$ to $A$
has at least two global minima.}\par
\smallskip
PROOF. For each $(x,\lambda)\in X\times {\Lambda}$, put
$$f(x,\lambda)=J(x)-\mu\varphi(\Psi(x,\lambda))\ .$$ 
By Theorem 4.1, there exists $A\in {\cal N}$ such that
$$\sup_{\lambda\in D}\inf_{x\in A}f(x,\lambda)<
\inf_{x\in A}\sup_{\lambda\in D}f(x,\lambda)\ ,$$
where
$$D=\hbox {\rm conv}(\{\lambda_x : x\in A\})\ .$$
Now, the conclusion comes directly applying Theorem 1.2 to
the restriction of $f$ to $A\times D$.\hfill $\bigtriangleup$\par
\medskip
A non-empty set $C$ in a normed space $S$ is said to be
uniquely remotal with respect to a set $D\subseteq S$ if, for each $y\in D$, there exists a unique $x\in C$
such that
$$\|x-y\|=\sup_{u\in C}\|u-y\|\ .$$
The main problem in theory of such sets is to know if they are singletons.\par
\smallskip
If $E, F$ are two real vector spaces and $D$ is a convex subset of $E$, 
we say that an operator $\Phi:D\to F$ is affine if
$$\Phi(\lambda x+(1-\lambda)y)=\lambda \Phi(x)+(1-\lambda)\Phi(y)$$
for all $x, y\in D$, $\lambda\in [0,1]$.
\smallskip
The next three results are applications of Theorem 5.4.\par
\medskip
THEOREM 5.5. - {\it Let $Y$ be a real normed space and let $X\subseteq Y$
be a non-empty compact uniquely
remotal set with respect to \hbox {\rm conv}$(X)$.\par
Then, $X$ is a singleton.}\par
\smallskip
PROOF. Arguing by contradiction, assume that $X$ contains at least 
two points.
Now, apply Theorem 5.4 taking: $\Lambda=Y$, $y_0=0$, $\varphi(x)=\|x\|$, $\Psi(x,\lambda)=
x-\lambda$, $J=0$ and ${\cal N}=\{X\}$. Note that we are allowed to apply Theorem 5.4 since
$x\to \lambda_x$ is not constant. Then, it would exists
 $\lambda^*\in \hbox {\rm conv}(X)$
such that the function $x\to -\|x-\lambda^*\|$ has at least two global minima in
$X$, against the hypotheses.\hfill $\bigtriangleup$\par
\medskip
REMARK 5.2. - Observe that Theorem 5.5 improves a classical result by
V. L. Klee ([5]) under
two aspects: $Y$ does not need to be complete and $\overline {\hbox {\rm conv}}(X)$
is replaced by $\hbox {\rm conv}(X)$. Note also that our proof is completely different
from that of Klee which is based on the Schauder fixed point theorem.
\medskip
THEOREM 5.6. - {\it Let $X$ be a finite-dimensional real Hilbert space and 
$J:X\to {\bf R}$ a $C^1$ function. Set
$$\eta=\liminf_{\|x\|\to +\infty}{{J(x)}\over {\|x\|^2}}$$
and
$$\theta=\inf\left \{
{{J(x)-J(u)}\over {\|x-u\|^2}} : (u,x)\in M_{J}\times
X\hskip 3pt \hbox {\rm with}\hskip 3pt x\neq u\right
\}$$
where $M_J$ denotes the set of all global minima of $J$.
Assume that $$\theta<\eta\ .$$
Then, for each $\mu\in ]2\theta,2\eta[$, there exists $y_{\mu}\in X$ such that
the equation
$$J'(x)-\mu x=y_{\mu}$$
has at least three solutions.}\par
\smallskip
PROOF. 
Let $\mu\in ]2\theta,2\eta[$. We clearly
have
$$\lim_{\|x\|\to +\infty}\left ( J(x)-{{\mu}\over {2}}\|x-
\lambda\|^2\right ) =+\infty
\eqno{(5.3)}$$
for all $\lambda\in X$. So, since $X$ is finite-dimensional,
the function $x\to
J(x)-{{\mu}\over {2}}\|x-\lambda\|^2$ is continuous and inf-compact for
all $\lambda\in X$. Therefore, we can apply Theorem 5.4 
 taking: $X=Y=\Lambda$, $y_0=0$, $\varphi(y)=\|y\|^2$,
$\Psi(x,\lambda)=x-\lambda$ and ${\cal N}=\{X\}$. Consequently, there exists $\lambda^*_{\mu}\in
X$, such that the function $x\to J(x)-{{\mu}\over
{2}}\|x-\lambda^*_{\mu}\|^2$ has at least two global minima. 
By $(5.3)$ and the finite-dimensionality of $X$ again, the same function
satisfies the Palais-Smale condition, and so it admits at least three critical
points, thanks to Corollary 1 of [9]. Of course, this gives the conclusion, taking
$y_{\mu}=\lambda^*_{\mu}$.
\hfill $\bigtriangleup$\par
\medskip
REMARK 5.3. - Clearly, there are two situations in which Theorem 5.6 can immediately
be applied: when $\eta=+\infty$, and when $\eta>0$ and $\theta=0$. Note
that one has $\theta=0$ if, in particular, $J$ possesses at least two global
minima. \par
\medskip
REMARK 5.4. - It is also clear that under the same assumptions as those of Theorem 5.6
 but the finite-dimensionality of $X$, the conclusion is still true 
for every $\mu\in ]2\theta,2\eta[$ such that, for each $\lambda\in X$, the functional
$x\to J(x)-{{\mu}\over {2}}\|x-\lambda\|^2$ is weakly lower
semicontinuous and satisfies the Palais-Smale condition.\par
\medskip
THEOREM 5.7. - {\it Let $Y$ be a 
finite-dimensional real Hilbert space,
$J:Y\to {\bf R}$ a $C^1$ function with locally Lipschitzian derivative,
and $\varphi:Y\to [0,+\infty[$ a $C^1$ convex function with locally
lipschitzian derivative at $0$ and $\varphi^{-1}(0)=\{0\}$.\par
Then, for each $x_0\in Y$ for wich $J'(x_0)\neq 0$, there exists $\delta>0$
such that, for each $r\in ]0,\delta[$, the restriction of $J$ to
$B(x_0,r)$ has a unique global minimum $u_r$ which
satisfies
$$J(u_r)\leq J(x)-\varphi(x-u_r)$$
for all $x\in B(x_0,r)$,
where
$$B(x_0,r)=\{x\in Y: \|x-x_0\|\leq r\}\ .$$}
\medskip
PROOF. First of all, observe that $\varphi\in {\cal
G}$, with $y_0=0$. Indeed, let $V\subset Y$ be a neighbourhood of
$0$ and let $s>0$ be such that $B(0,s)\subseteq V$. Set
$$\alpha=\inf_{\|x\|=s}\varphi(x)\ .$$
Since dim$(Y)<\infty$, $\partial B(0,s)$ is compact and so $\alpha>0$.
Let $x\in Y$ with $\|x\|>s$. Let $S$ be the segment joining $0$ and $x$.
By convexity, we have
$$\varphi(z)\leq \varphi(x)$$
for all $z\in S$. Since $S$ meets $\partial B(0,s)$, we infer that
$\alpha\leq \varphi(x)$. Hence, we have
$$\alpha\leq \inf_{\|x\|>s}\varphi(x)\leq \inf_{x\in X\setminus V}\varphi(x)
\ .$$
Now, fix $x_0\in Y$ with $J'(x_0)\neq 0$. Taking into account that
$\varphi'(0)=0$, by continuity, we can choose $\sigma>0$ so that
$$\|\varphi'(\lambda)\|<\|J'(x)\|$$
 for all $(x,\lambda)\in B(x_0,\sigma)\times B(0,\sigma)$.
For each $(x,\lambda)\in Y\times Y$,
put
$$f(x,\lambda)=J(x)-\varphi(x-x_0-\lambda)\ .$$
Of course, we have
$$f'_x(x,\lambda)\neq 0$$
for all $(x,\lambda)\in B\left ( x_0,{{\sigma}\over {2}}\right ) 
\times B\left ( 0,{{\sigma}\over {2}}\right )$.
Next, since $J'$ is locally Lipschitzian at $x_0$  and $\varphi'$ is
locally Lipschitzian at $0$, there are $\rho\in \left 
] 0,{{\sigma}\over {2}}\right ] $ and
$L>0$ such that
$$\|f'_x(x,\lambda)-f'_x(y,\lambda)\|\leq L\|x-y\|$$
for all $x,y\in B(x_0,\rho)$, $\lambda\in B(0,\rho)$. Now,
fix $\lambda\in B(0,\rho)$. Denote by $\Gamma_{\lambda}$ the set
of all global minima of the restriction of the
function $x\to f(x,\lambda)+{{L}\over {2}}\|x-x_0\|^2$ to 
$B(x_0,\rho)$. Note that $x_0\not\in \Gamma_{\lambda}$ (since
$f'_x(x_0,\lambda)\neq 0$). As $f$ is
continuous, the multifunction $\lambda\to \Gamma_{\lambda}$ is
upper semicontinuous and so the function $\lambda\to \hbox
{\rm dist}(x_0,\Gamma_{\lambda})$ is lower semicontinuous. As
a consequence, by compactness, we have
$$\delta:=\inf_{\lambda\in B(0,\rho)}
\hbox {\rm dist}(x_0,\Gamma_{\lambda})>0\ .$$
At this point, from the proof of Theorem 1 of [11] it follows that, for
each $\lambda\in B(0,\rho)$ and each $r\in ]0,\delta[$, 
the restriction of function $f(\cdot,\lambda)$ to $B(x_0,r)$ has a
unique global minimum. Fix $r\in ]0,\delta[$. Apply Theorem 5.4
with $X=B(x_0,r)$, $\Lambda=Y$, ${\cal N}=\{B(x_0,r)\}$ and $\Psi(x,\lambda)=
x-x_0-\lambda$. With such choices, its conclusion does not hold with
$\mu=1$ (recall, in particular, that
$r<\rho$). This implies that $1\leq \theta(\varphi,\Psi,J)$ since the other
assumptions are satisfied. But the above inequality is just equivalent to
$$J(u_r)\leq J(x)-\varphi(x-u_r)$$
for all $x\in B(x_0,r)$, where $u_r$ is the unique global minimum of
$J_{|B(x_0,r)}$, and the proof is complete.\hfill $\bigtriangleup$\par
\medskip
A joint application of Theorems 1.1 and 4.1 gives\par
\medskip
THEOREM 5.8. - {\it Let $\varphi\in {\cal G}$, $\Psi\in {\cal H}$  
 and $J\in {\cal M}$. 
Moreover, assume that $X$ is a topological space, that
$\Lambda$ is a real topological vector space and that
$\varphi(\Psi(x,\cdot))$ is quasi-convex and continuous for each $x\in X$. Finally, let
$\mu>\theta(\varphi,\Psi,J)$  and let $C\subseteq \Lambda$ be a convex set, with
$\{\lambda_x : x\in X\}\subseteq \overline {C}$,
such that
 the function $x\to J(x)-\mu\varphi(\Psi(x,\lambda))$ is
lower semicontinuous and inf-compact in $X$ for all $\lambda\in C$.\par
Under such hypotheses, 
there exists $\lambda^*\in C$
 such
that the function $x\to J(x)-\mu\varphi(\Psi(x,\lambda^*))$ 
has at least two global minima in $X$.}\par
\smallskip
PROOF. Set
$$D=\{\lambda_x : x\in X\}$$
and, for each $(x,\lambda)\in X\times \Lambda$, put
$$f(x,\lambda)=J(x)-\mu\varphi(\Psi(x,\lambda))\ .$$
Theorem 4.1 ensures that
$$\sup_{\lambda\in \Lambda}\inf_{x\in X}f(x,\lambda)<\inf_{x\in X}\sup_{\lambda\in D}f(x,\lambda)\ .\eqno{(5.4)}$$
But, since $f(x,\cdot)$ is continuous and $D\subseteq \overline {C}$,
 we have
$$\sup_{\lambda\in D}f(x,\lambda)=\sup_{\lambda\in \overline {D}}f(x,\lambda)\leq \sup_{\lambda\in \overline {C}}f(x,\lambda)
= \sup_{\lambda\in C}f(x,\lambda)$$
for all $x\in X$, and hence, from $(5.4)$, it follows that 
$$\sup_{\lambda\in C}\inf_{x\in X}f(x,\lambda)<
\inf_{x\in X}\sup_{\lambda\in D}f(x,\lambda)\leq \inf_{x\in X}\sup_{\lambda\in C}f(x,\lambda)\ .$$
At this point, the conclusion follows applying Theorem 1.1 to the restriction of
the function $f$ to $X\times C$.\hfill $\bigtriangleup$\par
\medskip
The next result comes from a joint application of Theorems 4.3 and 1.C.\par
\medskip
THEOREM 5.9. -  {\it Let $X$ be a real inner product space and let $\tau$ be a topology on $X$.
Moreover, let $J:X\to {\bf R}$ be a functional such that
$$J(0)=0<\sup_X J$$
and
$$\beta^*:=\sup_{x\in X\setminus \{0\}}{{J(x)}\over {\|x\|^2}}<+\infty\ .\eqno{(5.5)}$$
Finally, let $\lambda>{{1}\over {\beta^*}}$ and let ${\cal N}$ be a weakly filtering cover of $X$ such that, for
each $A\in {\cal N}$ and each $y\in X$, the restriction to $A$ of the functional $x\to \|x\|^2-\lambda J(x)+\langle x,y\rangle$ is
$\tau$-lower semicontinuous and inf-$\tau$-compact.\par
Then, there exists $\tilde A\in {\cal N}$ with the following property: for every convex set $C\subseteq X$ whose closure (in the strong 
topology) contains $\tilde A$, there exists $\tilde y\in C$ such that the restriction to $\tilde A$ of the functional 
$x\to \|x\|^2-\lambda J(x)+\left \langle x,2(\beta^*\lambda-1)\tilde y\right \rangle$ has
at least two global minima.}\par
\smallskip
PROOF. In view of $(5.5)$, we have 
$$\inf_{x\in X}(\|x\|^2-\lambda J(x))<0 \ ,$$
as well as
$$\inf_{x\in X}(\|x\|^2-\lambda J(x)+(\beta^*\lambda-1)\|x\|^2)\geq 0\ .$$
So, we can apply Theorem 4.3 taking $Y=X$,
$$\mu=\beta^*\lambda-1\ ,$$
$$I(x)=\|x\|^2-\lambda J(x)$$
and 
$$\Phi(x)=x\ .$$
Therefore, there exists $\tilde A\in {\cal N}$ such that
$$\sup_{y\in Y}\inf_{x\in \tilde A}(\|x\|^2-\lambda J(x)+(\beta^*\lambda-1)(2\langle x,y\rangle-
\|y\|^2))<
\inf_{x\in \tilde A}\sup_{y\in \tilde A}(\|x\|^2-\lambda J(x)+(\beta^*\lambda-1)(2\langle x,y\rangle-\|y\|^2))\ .\eqno{(5.6)}$$\par
Now, consider the function $f:X\times X\to {\bf R}$ defined by
$$f(x,y)=\|x\|^2-\lambda J(x)+(\beta^*\lambda-1)(2\langle x,y\rangle-\|y\|^2)$$
for alla $(x,y)\in X\times X$. 
 Since $f(x,\cdot)$ is continuous and $\tilde A\subseteq \overline {C}$,
 we have
$$\sup_{y\in \tilde A}f(x,y)=\sup_{v\in \overline {\tilde A}}f(x,y)\leq \sup_{v\in \overline {C}}f(x,y)
= \sup_{y\in C}f(x,y)$$
for all $x\in X$, and hence, taking $(5.6)$ into account, it follows that 
$$\sup_{y\in C}\inf_{x\in \tilde A}f(x,y)<
\inf_{x\in \tilde A}\sup_{y\in \tilde A}f(x,y)\leq \inf_{x\in \tilde A}\sup_{y\in C}f(x,y)\ .\eqno{(5.7)}$$
Now, in view of $(5.7)$, taking into account that $f_{|\tilde A\times C}$ is $\tau$-lower semicontinuous and inf-$\tau$-compact in $\tilde A$, and 
continuous
and concave in $C$, we can apply Theorem 1.C to $f_{|\tilde A\times C}$. Consequently, there exists $\tilde y\in C$ such that
$f_{|\tilde A}(\cdot,\tilde y)$ has at least two global minima, and the proof is complete. \hfill $\bigtriangleup$\par
\medskip
In turn, from Theorem 5.9, we get\par
\medskip
THEOREM 5.10. - {\it Let $X$ be a real Hilbert space and let $J:X\to {\bf R}$ be a $C^1$ functional, with compact derivative,
such that 
$$\alpha^*:=\max\left \{0,\limsup_{\|x\|\to +\infty}{{J(x)}\over {\|x\|^2}}\right \}<\beta^*:=\sup_{x\in X\setminus \{0\}}{{J(x)}\over 
{\|x\|^2}}<+\infty\ .$$
Then, for every $\lambda\in \left ]{{1}\over {2\beta^*}}, {{1}\over {2\alpha^*}}\right [$ and
for every convex set $C\subseteq X$ dense in $X$, there exists $\tilde y\in C$ such that the equation
$$x=\lambda J'(x)+\tilde y$$
has at least three solutions, two of which are global minima of the functional
$x\to {{1}\over {2}}\|x\|^2-\lambda J(x)-\langle x,\tilde y\rangle$\ .}\par
\smallskip
PROOF. Fix $\lambda\in \left ]{{1}\over {2\beta^*}}, {{1}\over {2\alpha^*}}\right [$ and
a convex set $C\subseteq X$ dense in $X$. For each $y\in X$, we have
$$\liminf_{\|x\|\to +\infty}\left ( 1-2\lambda{{J(x)}\over {\|x\|^2}}-{{\langle x,y\rangle}\over {\|x\|^2}}\right )=
1-2\lambda\limsup_{\|x\|\to +\infty}{{J(x)}\over {\|x\|^2}}>0\ .$$
So, from the identity
$$\|x\|^2-2\lambda J(x)-\langle x,y\rangle=\|x\|^2\left ( 1-2\lambda{{J(x)}\over {\|x\|^2}}-{{\langle x,y\rangle}\over {\|x\|^2}}\right )$$
it follows that
$$\lim_{\|x\|\to +\infty} (\|x\|^2-2\lambda J(x)-\langle x,y\rangle)=+\infty\ .\eqno{(5.8)}$$
Since $J'$ is compact, $J$ is sequentially weakly continuous ([29], Corollary 41.9).  Then, in view of $(5.8)$ and of the Eberlein-Smulyan theorem, 
for each
$y\in X$, the functional $x\to \|x\|^2-2\lambda J(x)+\langle x,y\rangle$ is inf-weakly compact in $X$. So, we can apply Theorem 5.9 taking the weak 
topology as $\tau$
and ${\cal N}=\{X\}$. Consequently, since the set ${{1}\over {1-2\beta^*\lambda}}C$ is convex and dense in $X$, there exists 
$\hat y\in {{1}\over {1-2\beta^*\lambda}}C$ such that the functional  
$x\to \|x\|^2-2\lambda J(x)+\left \langle x,2(2\beta^*\lambda-1)\hat y\right \rangle$ has
at least two global minima in $X$ which are two of its critical points. Since the same functional satisfies the Palais-Smale 
condition ([29], Example 38.25), it has a third critical 
point in view of Corollary 1 of [9]. Clearly, the conclusion follows taking $\tilde y=(1-2\beta^*\lambda)\hat y$\ .\hfill $\bigtriangleup$\par
\medskip
Let us conclude this section with a further consequence of Theorem 1.2.\par
\smallskip
Let us introduce the following notations. We denote by ${\bf R}^X$ the space of
all functionals $\varphi:X\to {\bf R}$.  For each $I\in {\bf R}^X$ and for
each  of non-empty subset $A$ of $X$, we denote by
$E_{I,A}$ the set of all $\varphi\in {\bf R}^X$ such that $I+\varphi$ is sequentially
weakly lower semicontinuous and coercive, and
$$\inf_A\varphi\leq 0\ .$$
\medskip
THEOREM 5.11. - {\it Let  $I:X\to {\bf R}$ be a functional and $A, B$ two non-empty subsets of $X$ such that
$$\sup_AI<\inf_BI\ .\eqno{(5.9)}$$
Then, for every convex set $Y\subseteq E_{I,A}$ such that
$$\inf_{x\in B}\sup_{\varphi\in Y}\varphi(x)\geq 0\hskip 5pt and\hskip 5pt  
\inf_{x\in X\setminus B}\sup_{\varphi\in Y}\varphi(x)=+\infty\ ,\eqno{(5.10)}$$
 there exists $\tilde\varphi\in Y$ such that the
functional $I+\tilde\varphi$ has at least two global minima.}\par
\smallskip
PROOF. 
Consider the function $f:X\times {\bf R}^X\to {\bf R}$ defined by
$$f(x,\varphi)=I(x)+\varphi(x)$$
for all $x\in X$, $\varphi\in {\bf R}^X$. Fix $\varphi\in Y$. In view of $(5.9)$, we also can fix $\epsilon\in ]0,\inf_BI-\sup_AI[$.
Since $\inf_A\varphi\leq 0$, there is $\bar x\in A$ such that $\varphi(\bar x)<\epsilon$. Hence,
we have
$$\inf_{x\in X}(I(x)+\varphi(x))\leq I(\bar x)+\varphi(\bar x)<\sup_AI+\epsilon\ ,$$
from which it follows that
$$\sup_{\varphi\in Y}\inf_{x\in X}(I(x)+\varphi(x))\leq \sup_AI+\epsilon<\inf_BI\ .\eqno{(5.11)}$$
On the other hand, in view of $(5.10)$, one has
$$\inf_BI\leq \inf_{x\in B}(I(x)+\sup_{\varphi\in Y}\varphi(x))=\inf_{x\in B}\sup_{\varphi\in Y}(I(x)+\varphi(x))=
\inf_{x\in X}\sup_{\varphi\in Y}(I(x)+\varphi(x))\ .\eqno{(5.12)}$$
Finally, from $(5.11)$ and $(5.12)$, it follows that
$$\sup_{\varphi\in Y}\inf_{x\in X}f(x,\varphi)<\inf_{x\in X}\sup_{\varphi\in Y}f(x,\varphi)\ .$$
Therefore, the function $f$ satisfies the assumptions of Theorem 1.2, and the conclusion follows.\hfill $\bigtriangleup$
\medskip
Notice the following remarkable corollary of Theorem 5.11:\par
\medskip
COROLLARY 5.1. - {\it Let $I:X\to {\bf R}$ be a sequentially weakly lower semicontinuous, non-convex functional such that
$I+\varphi$ is coercive for all $\varphi\in X^*$.\par
Then, for every convex set $Y\subseteq X^*$, dense in $X^*$, there exists $\tilde\varphi\in Y$ such that the functional
$I+\tilde\varphi$ has at least two global minima.}\par
\smallskip
PROOF. Since $I$ is not convex,  there
exist $x_1, x_2\in X$ and $\lambda\in ]0,1[$ such that
$$\lambda I(x_1)+(1-\lambda) I(x_2)<I(x_3)$$
where
$$x_3=\lambda x_1+(1-\lambda)x_2\ .$$
Fix $\psi\in X^*$ so that
$$\psi(x_1)-\psi(x_2)=I(x_1)-I(x_2)$$
and put
$$\tilde I(x)=I(x_3-x)-\psi(x_3-x)$$
for all $x\in X$. It is easy to check that
$$\tilde I(\lambda(x_1-x_2))=\tilde I((1-\lambda)(x_2-x_1))<\tilde I(0)\ .\eqno{(5.13)}$$
Fix a
convex set $Y\subseteq X^*$ dense in $X^*$ and put
$$\tilde Y=-Y-\psi\ .$$
Hence, $\tilde Y$ is convex and dense in $X^*$ too.  Now, set
$$A=\{\lambda(x_1-x_2),(1-\lambda)(x_2-x_1)\}\ .$$
Clearly, we have
$$X^*\subset E_{\tilde I, A}\ .\eqno{(5.14)}$$ 
Since $\tilde Y$ is dense in $X^*$, we have
$$\sup_{\varphi\in \tilde Y}\varphi(x)=+\infty$$
for all $x\in X\setminus \{0\}$. Hence, in view of $(5.13)$ and $(5.14)$, we can apply
Theorem 5.11 with  $B=\{0\}$, $I=\tilde I$, $Y=\tilde Y$. Accordingly,
there exists $\tilde\varphi\in Y$ such that the functional $\tilde I-\tilde\varphi-\psi$
has two global minima in $X$, say $u_1, u_2$. At this point, it is clear that
$x_3-u_1, x_3-u_2$ are two global minima of the functional $I+\tilde\varphi$,
and the proof is complete.\hfill $\bigtriangleup$

\bigskip
{\bf 6. A range property for non-expansive potential operators}\par
\bigskip
In this section, $(X,\langle\cdot,\cdot\rangle)$ is an infinite-dimensional
real Hilbert space and $T:X\to X$ is a non-expansive potential operator. This means
that
$$\|T(x)-T(y)\|\leq \|x-y\|$$
for all $x,y\in X$ and that $T$ is the G\^ateaux derivative of a functional
$J:X\to {\bf R}$. \par
\smallskip
For instance, any continuous symmetric linear operator from $X$ into itself, with
norm less than or equal to $1$, is a non-expansive potential operator.\par
\smallskip
Another classical example of such operators is as follows. Let 
$f:X\to {\bf R}$ be a convex continuous function and, for each
$x\in X$, let 
$\partial f(x)$ denote the sub-differential of $f$ at $x$, i. e.
$$\partial f(x)=\left \{ z\in X :
\inf_{y\in X}(f(y)-\langle z,y\rangle)\geq f(x)-\langle z,x\rangle
\right \}\ .$$
Then $x\to x-(\hbox {\rm id}+\partial f)^{-1}(x)$ is a non-expansive
potential operator.\par
\smallskip
Now let $\Phi:X\to X$ be the operator defined by
$$\Phi(x)=x+T(x)$$
for all $x\in X$.\par
\smallskip
The following result does highlight a range property of the operator $\Phi$. The proof is
based on combining some ideas from [11] with Theorem 1.C.
\par
\medskip
THEOREM 6.1. - {\it If the functional $J$ is sequentially weakly lower
semicontinuous, then 
there exists a closed ball $B$ in $X$ such that
$\Phi(B)$ intersects each convex and dense subset of $X$\ .}\par
\medskip
Before proving Theorem 6.1, some remarks are in order.\par
\smallskip
When $T$ is a contraction, Theorem 6.1 is immediate. Actually, in that case,
thanks to the Banach fixed point principle, the operator $\Phi$ turns out to be
a homeomorphism between $X$ and itself. Hence, the really interesting case is when the Lipschitz
constant of $T$ is exactly $1$.\par
\smallskip
Theorem 6.1 is no longer true if $J$ is not sequentially weakly lower semicontinuous.
In this connection, the simplest example is provided by $T(x)=-x$. Actually,
since dim$(X)=\infty$, the norm is not sequentially weakly upper
semicontinuous.
\smallskip
A further remark is that, under the assumptions of Theorem 6.1, it may happen that
the set $\Phi(B)$ has an empty interior for every ball $B$ in $X$. In this
connection,
consider the case where $T$ is a compact, symmetric, negative linear operator
with norm $1$.
In such a case, by classical results, $J$ is sequentially weakly
continuous and $\Phi(X)\neq X$. Since $\Phi(X)$ is a linear subspace, this
clearly implies that int$(\Phi(X))=\emptyset$.\par
\medskip
{\it Proof of Theorem 6.1.}
If the functional $J$ is convex, then $T$ is maximal monotone and, by a
classical result of Minty, $\Phi$ turns out to be a homeomorphism between
$X$ and itself. So, in that case, we are done.
Therefore, assume that $J$ is not convex. As a consequence, there
exist $x_1, x_2\in X$ and $\lambda\in ]0,1[$ such that
$$\lambda J(x_1)+(1-\lambda) J(x_2)<J(x_3)$$
where
$$x_3=\lambda x_1+(1-\lambda)x_2\ .$$
Fix $z\in X$ so that
$$\langle x_1-x_2,z\rangle=J(x_1)-J(x_2)$$
and put
$$\tilde J(x)=J(x_3-x)-\langle x_3-x,z\rangle$$
for all $x\in X$. Note that
$$\tilde J(\lambda(x_1-x_2))=\tilde J((1-\lambda)(x_2-x_1))<\tilde J(0)\ .\eqno{(6.1)}$$
Now, put
$$r=\max\{\lambda, 1-\lambda\}\|x_1-x_2\|$$
and denote by $C$ the closed ball in $X$ of radius $r$ centered at $0$.
Fix a
convex and dense set $V\subseteq X$ and put
$$Y=V-x_3-z\ .$$
Hence, $Y$ is convex and dense too.
Consider
the function $f:X\times Y\to {\bf R}$ defined by
$$f(x,y)=\tilde J(x)+\langle x,y\rangle$$
for all $(x,y)\in X\times Y$. Observing that, for each $y\in Y$, one has
 $$\min\{\langle \lambda(x_1-x_2),y\rangle, \langle (1-\lambda)(x_2-x_1),
y\rangle\}\leq 0\ ,$$ in view of the equality in $(6.1)$, it follows that
$$\sup_{y\in Y}\inf_{x\in C}f(x,y)\leq
\tilde J(\lambda(x_1-x_2))\ .\eqno{(6.2)}$$
On the other hand, by the density of $Y$, for each $x\in C\setminus
\{0\}$, one has
$$\sup_{y\in Y}\langle x,y\rangle=+\infty$$
and hence
$$\inf_{x\in C}\sup_{y\in Y}f(x,y)=\tilde J(0)\ .\eqno{(6.3)}$$
Thus, from $(6.1), (6.2), (6.3)$ it follows that
$$\sup_{y\in Y}\inf_{x\in C}f(x,y)<
\inf_{x\in C}\sup_{y\in Y}f(x,y)\ .$$
Then, since $f(\cdot,y)_{|C}$ is weakly
lower semicontinuous in $C$ (thanks to the Eberlein-Smulyan
theorem) and $f(x,\cdot)$ is concave
and continuous in $Y$, we can apply Theorem 1.C. Accordingly,
there exists $\hat y\in Y$ such that $f(\cdot,\hat y)_{|C}$
has at least two global minima $u_1, u_2$ in $C$. Now,
consider the function $g:X\times [1,+\infty[\to {\bf R}$ defined by
$$g(x,\lambda)=
{{\lambda}\over {2}}(\|x\|^2-r^2)+f(x,\hat y)$$
for all $(x,\lambda)\in X\times [1,+\infty[$. Observe that
 the functional $g(\cdot,\lambda)$ (besides being
 continuous) is strictly convex and coercive
if $\lambda>1$, while it is convex if $\lambda=1$. Indeed, let
$\lambda\geq 1$.
For each $x, y\in X$, we have\par
$$\langle \lambda x+T(x_3-x)-\lambda y+T(x_3-y),x-y\rangle=
\lambda\|x-y\|^2-\langle
T(x_3-x)-T(x_3-y),x-y\rangle\geq$$
$$\lambda\|x-y\|^2-\|T(x_3-x)-T(x_3-y)\|\|x-y\|
\geq (\lambda-1)\|x-y\|^2\ .$$
 From this, it follows that the G\^ateaux derivative of the functional
$g(\cdot,\lambda)$
 (that is, the operator $x\to \lambda x-T(x_3-x)+z+\hat y$) 
 is monotone
and that it is uniformly monotone if $\lambda>1$. Now the claim
follows from classical results ([29], pp. 247-248). Furthermore, for
each $x\in X$, the function $g(x,\cdot)$ is concave and continuous, and
$\lim_{\lambda\to +\infty}g(0,\lambda)=-\infty$. So, we are allowed to
apply a classical saddle-point theorem ([29], Theorem 49.A) to
the function $g$. Accordingly,
there exists $(\hat x,\hat \lambda)\in X\times [1,+\infty[$ such
that
$$g(\hat x,\hat\lambda)=\inf_{x\in X}g(x,\hat\lambda)=\sup_{\lambda\geq 1}
g(\hat x,\lambda)\ .$$
This implies that $\hat x\in C$,
$${{\hat\lambda}\over {2}}\|\hat x\|^2+
f(\hat x,\hat y)=
\inf_{x\in X}\left ( 
{{\hat\lambda}\over {2}}\|x\|^2+f(x,\hat y)\right )\eqno{(6.4)}$$
and
$${{\hat\lambda}\over {2}}(\|\hat x\|^2-r^2)=
{{1}\over {2}}(\|\hat x\|^2-r^2)\ .\eqno{(6.5)}$$
We claim that $\hat\lambda=1$. If $\|\hat x\|<r$, this follows directly
from $(6.5)$. So, assume that $\|\hat x\|=r$. In this case, for $i=1, 2$, we
have
$${{\hat \lambda}\over {2}}\|u_i\|^2+f(u_i,\hat y)\leq
{{\hat \lambda}\over {2}}\|\hat x\|^2+f(\hat x,\hat y)$$
and hence from $(6.4)$ it follows that
$${{\hat \lambda}\over {2}}\|u_i\|^2+f(u_i,\hat y)=
\inf_{x\in X}\left ( 
{{\hat\lambda}\over {2}}\|x\|^2+f(x,\hat y)\right )\ .\eqno{(6.6)}$$
But, for $\lambda>1$, the functional $x\to {{\lambda}\over {2}}\|x\|^2
+f(x,\hat y)$ has a unique global minimum in $X$ because it is strictly convex.
So, the equality $\hat \lambda=1$ follows from $(6.6)$. Therefore, by $(6.4)$,
the G\^ateaux derivative of the functional $x\to {{1}\over {2}}\|x\|^2+
f(x,\hat y)$ vanishes at $\hat x$. This means that
$$T(x_3-\hat x)-\hat x=z+\hat y\ .$$
Therefore, if $B$ is the closed ball of radius $r$ centered at $x_3$, we have
$x_3-\hat x\in B$ and $\Phi(x_3-\hat x)\in V$, and the proof is complete.
\hfill $\bigtriangleup$
\bigskip
{\bf 7. Singular points of non-monotone potential operators}\par
\bigskip
In this section, $(X,\|\cdot\|)$ is a reflexive real
Banach space, with topological dual $X^*$, and $T:X\to X^*$ is a continuous potential operator. 
As a consequence, the functional
$$x\to J_T(x):=\int_0^1T(sx)(x)ds$$
is of class $C^1$ in $X$ and its G\^ateaux derivative is equal to $T$.
\smallskip
Let us recall a few classical definitions.\par 
\smallskip
$T$ is said to be monotone if
$$(T(x)-T(y))(x-y)\geq 0$$
for all $x, y\in X$.  This is equivalent to the fact that the functional $J_T$ is convex.
\smallskip
$T$ is said to be closed if for each closed set $C\subseteq X$, the set $T(C)$ is closed in $X^*$.\par
\smallskip
$T$ is said to be compact if for each bounded set $B\subset X$, the set $\overline {T(B)}$ is compact in $X^*$.\par
\smallskip
$T$ is said to be proper if for each compact set $K\subset X^*$, the set $T^{-1}(K)$ is compact in $X$.\par
\smallskip
$T$ is said to be a local homeomorphism at
a point $x_0\in X$ if there are a neighbourhood $U$ of $x_0$ and a neighbourhood
$V$ of $T(x_0)$ such that the restriction of $T$ to $U$ is a homeomorphism between $U$
and $V$. If $T$ is not a local homeomorphism at $x_0$, we say that $x_0$ is
a singular
point of $T$.\par
\smallskip
We denote by $S_T$ the set of all singular points of $T$. Clearly, the set $T$ is closed.\par
\smallskip
Assume that the restriction of $T$ to some open set $A\subseteq X$ is of
class $C^1$.\par
\smallskip
 We  then denote by
 $\tilde S_{T_{|A}}$ the set of all $x_0\in A$ such that
the operator
$T'(x_0)$ is not invertible. Since the set of all invertible operators belonging to ${\cal L}(X,X^*)$ is
open in ${\cal L}(X,X^*)$, by the continuity of $T'$,
the set $\tilde S_{T_{|A}}$ is closed in $A$.\par
\smallskip
Also,  $T$ is said to be a Fredholm operator of index zero in $A$
if, for each $x\in A$, the codimension of $T'(x)(X)$ and  the dimension of
$(T'(x))^{-1}(0)$ are finite and equal.
\smallskip
A set in a topological space is said to be $\sigma$-compact if it is the union of an at most countable
family of compact sets.\par
\smallskip
A functional $I:X\to {\bf R}$ is said to be coercive if
$$\lim_{\|x\|\to +\infty}I(x)=+\infty\ .$$
Let us recall the two following results:\par
\medskip
THEOREM 7.A. - ([23], Theorem 2.1). - {\it If $X$ is infinite-dimensional, if  $T$
is closed and if $S_T$ is $\sigma$-compact, then the
restriction of $T$ to $X\setminus S_T$ is a homeomorphism
between $X\setminus S_T$ and $X\setminus T(S_T)$.}\par
\medskip
THEOREM 7.B. - ([8], Theorem 5). - {\it If \hbox {\rm dim}$(X)\geq 3$,
if $T$ is a $C^1$
proper Fredholm operator of index zero and if $\tilde S_T$ is discrete,
then $T$ is a homeomorphism between $X$ and $X^*$.}\par
\medskip
We wish to show that a joint application of these results with Corollary 5.1 gives the following ones:
\medskip
THEOREM 7.1. - {\it If $X$ is infinite-dimensional, if $T$ is closed and non-monotone, if $J_T$ is
 sequentially weakly lower
semicontinuous and $J_T+\varphi$ is coercive
for all $\varphi\in X^{*}$, 
then 
both $S_T$ and $T(S_T)$ are not $\sigma$-compact.}\par
\medskip
THEOREM 7.2. - {\it In addition to the assumptions of Theorem 1, suppose that
there exists a
closed, $\sigma$-compact set $B\subset X$ such that the restriction
of $T$ to $X\setminus B$ is of class $C^1$.\par
Then, both $\tilde S_{T_{|(X\setminus B)}}$ and
$T(\tilde S_{T_{|(X\setminus B)}})$ are not $\sigma$-compact.}\par
\medskip
THEOREM 7.3. - {\it Assume that $(X,\langle\cdot,\cdot\rangle)$ is a Hilbert space, with $\hbox {\rm dim}(X)\geq 3$, and that $T$ is compact and
of class $C^1$ with
$$\liminf_{\|x\|\to +\infty}{{J_T(x)}\over {\|x\|^2}}\geq 0 \eqno{(7.1)}$$ 
and, for some $\lambda_0\geq 0$,
$$\lim_{\|x\|\to +\infty}\|x+\lambda T(x)\|=+\infty \eqno{(7.2)}$$
for all $\lambda>\lambda_0$ .\par
Set
$$\Gamma=\{(x,y)\in X\times X : \langle T'(x)(y), y\rangle<0\}$$
and, for each $\mu\in {\bf R}$, 
$$A_{\mu}=\{x\in X :  T'(x)(y)=\mu y\hskip 5pt for\hskip 5pt some\hskip 5pt y\in X\setminus \{0\}\}\ .$$
When $\Gamma\neq\emptyset$, set also
$$\tilde\mu=\max\left \{ -{{1}\over {\lambda_0}}, \inf_{(x,y)\in \Gamma}{{\langle T'(x)(y), y\rangle}\over {\|y\|^2}}\right
\}\ .$$
Then, the following assertions are equivalent:\par
\noindent
$(i)$\hskip 5pt the operator $T$ is not monotone\ ;\par
\noindent
$(ii)$\hskip 5pt there exists $\mu<0$ such that $A_{\mu}\neq\emptyset$\ ;\par
\noindent
$(iii)$\hskip 5pt $\Gamma\neq\emptyset$ and, for each $\mu\in ]\tilde\mu,0[$, the set
$A_{\mu}$ contains an accumulation point\ .}\par
\medskip
REMARK 7.1. - Of course, Theorem 7. 2 is meaningful only when $X$ and $X^*$ are linearly isomorphic. Indeed, if not, 
the fact that $\tilde S_{T_{|(X\setminus B)}}$  is not $\sigma$-compact follows directly from the equality
 $\tilde S_{T_{|(X\setminus B)}}=X\setminus B$ \ .\par
\medskip
We now establish the following technical proposition:\par
\medskip
PROPOSITION 7.1. - {\it If  $V$ is an infinite-dimensional real Banach space
space and if $U\subset V$ is a $\sigma$-compact set, then there exists
a convex cone $C\subset V$, dense in $V$, such that $U\cap C=\emptyset$.}\par
\smallskip
PROOF. Let us distinguish two cases.
First, assume that $V$ is separable. Fix a countable base $\{A_n\}$ of open sets in $V$. We claim
that there exists a sequence $\{x_n\}$ in $X$ such that, for each $n\in {\bf N}$,
 $$x_n\in A_n$$ and
$$U\cap C_{(x_1,...,x_n)}=\emptyset$$
where
$$C_{(x_1,...,x_n)}=\left \{ \sum_{i=1}^{n}\lambda_i x_i : 
\lambda_i \geq 0, \sum_{i=1}^{n}\lambda_i >0\right \} \ .$$
We proceed by induction on $n$. Clearly, the set $\cup_{\lambda>0}\lambda U$ is $\sigma$-compact
and so, since $X$ is infinite-dimensional, it does not contain 1.A. Thus, if we take $x_1\in
A_1\setminus \cup_{\lambda>0}\lambda U$, we have $U\cap C_{(x_1)} =\emptyset$.
Now, assume that $x_1,...,x_n$, with the desired properties, have been constructed. Consider
the set $\cup_{\mu>0}\mu(U-\overline {C_{(x_1,...,x_n)}})$. One readily sees that it is
$\sigma$-compact,
and so it does not
contain $A_{n+1}$. Choose
$x_{n+1}\in A_{n+1}\setminus \cup_{\mu>0}\mu(U-\overline {C_{(x_1,...,x_n)}})$.
Then, one has
$$U\cap C_{(x_1,...,x_{n+1})}=\emptyset\ .$$
Indeed, if there was $\hat x\in U\cap C_{(x_1,...,x_{n+1})}$, we would have
$\hat x=\sum_{i=1}^{n+1} \lambda_i x_i$,
with $\lambda_i\geq 0$ and $\sum_{i=1}^{n+1}\lambda_i>0$. In particular, $\lambda_{n+1}>0$, since
$U\cap C_{(x_1,...,x_n)}=\emptyset$. Consequently, we would have
$$x_{n+1}={{1}\over {\lambda_{n+1}}}\left ( \hat x-\sum_{i=1}^n \lambda_i x_i\right )$$
and so $x_{n+1}\in \cup_{\mu>0}\mu(U-\overline {C_{(x_1,...,x_n)}})$, against our choice. Thus,
the claimed sequence $\{x_n\}$ does exist. Now, put
$$C=\bigcup_{n=1}^{\infty}C_{(x_1,...,x_n)}\ .$$
It is clear that $C$ is a convex cone which does not meet $U$. Moreover, $C$ is dense in
$V$ since it meets each set $A_n$. Now, assume that $V$ is not separable. Let $\{x_{\gamma}\}_{\gamma\in \Gamma}$ 
be a Hamel basis of $V$. Set
$$\Lambda=\{\gamma\in \Gamma : x_{\gamma}\not\in \hbox {\rm span}(U)\}$$
and
$$L=\hbox {\rm span}(\{x_{\gamma} : \gamma\in \Lambda\})\ .$$
Clearly, span$(U)$ is separable since $U$ is so. Hence, $\Lambda$ is infinite.
Introduce in $\Lambda$ a total order $\leq$ with no greatest element. Next, for each $\gamma\in \Lambda$, let
$\psi_{\gamma}:L\to {\bf R}$ be a linear functional such that 
$$\psi_{\gamma}(x_{\alpha})=\cases {1 & if $\gamma=\alpha$\cr & \cr 0 & if $\gamma\neq\alpha$\ .\cr}$$
Now, set
$$D=\{x\in L : \exists \beta\in \Lambda : \psi_{\beta}(x)>0\hskip 5pt 
\hbox {\rm and}\hskip 5pt \psi_{\gamma}(x)=0\hskip 5pt \forall \gamma>\beta\}\ .$$
Of course, $D$ is a convex cone.
Fix $x\in L$. So, there is a finite set $I\subset\Lambda$ such that $x=\sum_{\gamma\in I}\psi_{\gamma}(x)x_{\gamma}$.
Now, fix $\beta\in\Lambda$ so that $\beta>\max I$. For each $n\in {\bf N}$, put
$$y_n=x+{{1}\over {n}}x_{\beta}\ .$$
Clearly, $\psi_{\beta}(y_n)={{1}\over {n}}$ and $\psi_{\gamma}(y_n)=0$ for all $\gamma>\beta$. Hence, $y_n\in D$.
Since $\lim_{n\to \infty}y_n=x$, we infer that $D$ is dense in $L$. At this point, it is immediate to check the
set $D+\hbox {\rm span}(U)$ is a convex cone, dense in $V$, which does not meet $U$.\hfill
$\bigtriangleup$\par
\medskip
{\it Proof of Theorem 7.1.} Let us prove that $S_T$ is not $\sigma$-compact.
 Arguing by contradiction, assume the contrary. Then, by Theorem 7.A,
for each $\varphi\in X^*\setminus T(S_T)$, the equation
$$T(x)=\varphi$$
has a unique solution in $X$.
Moreover, since $T$ is continuous, $T(S_T)$ is
 $\sigma$-compact too. Therefore, in view of Proposition 1, there is a convex set $Y\subset X^*$,
dense in $X^*$, such that $T(S_T)\cap Y=\emptyset$. On the other hand, thanks to Corollary 5.1, there is
$\tilde\varphi\in Y$ such that the functional $J_T-\tilde\varphi$ has at least two global minima in $X$ which
are therefore solutions of the equation
$$T(x)=\tilde\varphi\ ,$$
a contradiction. 
Now, let us prove that $T(S_T)$ is not $\sigma$-compact. Arguing by contradiction, assume the contrary.
 Consequently, since $T$ is proper ([23], Theorem 1),
 $T^{-1}(T(S_{T}))$ would
be $\sigma$-compact.  But then,  since $S_T$ is closed and $S_T\subseteq T^{-1}(T(S_T))$,  $S_T$ would be $\sigma$-compact,
a contradiction. The proof is complete.
 \hfill $\bigtriangleup$     \par
\medskip
{\it Proof of Theorem 7.2.}
By Theorem 7.1, the set $S_{T}$ is not $\sigma$-compact.
Now, observe that if
$x\in X\setminus {(\tilde S_{T_{|(X\setminus B)}}\cup B)}$, then,
by the inverse function
theorem, $T$ is a local homeomorphism at $x$, and so $x\not\in
S_{T}$. Hence, we have
$$S_{T}\subseteq \tilde S_{T_{|(X\setminus B)}}\cup B\ .$$
We then infer that $\tilde S_{T_{|(X\setminus B)}}$ is not
$\sigma$-compact since, otherwise,
$\tilde S_{T_{|(X\setminus B)}}\cup B$ would be so, and hence
also $S_{T}$ would be
$\sigma$-compact being closed. Finally, the fact that
$T(\tilde S_{T_{|(X\setminus B)}})$ is not $\sigma$-compact follows
as in the final part of the proof of Theorem 7.1, taking into account that
$\tilde S_{T_{|(X\setminus B)}}$ is closed in $X\setminus B$. \hfill $\bigtriangleup$\par
\medskip
{\it Proof of Theorem 7.3.} Clearly,  since $X$ is a Hilbert space, we are identifying $X^*$ to $X$.
 Let us prove that $(i)\to (iii)$. So, assume $(i)$.
Since $J_T$ is not convex, by a classical characterization ([27], Theorem 2.1.11), the set
$\Gamma$ is non-empty. Fix $\mu\in ]\tilde\mu,0[$. For each $x\in X$, put
$$I_{\mu}(x):={{1}\over {2}}\|x\|^2-{{1}\over {\mu}} J_T(x)\ .$$
Clearly, for some $(x,y)\in \Gamma$, we have
$$\left \langle y-{{1}\over {\mu}}T'(x)(y), y\right \rangle <0$$
and so, since 
$$I''_{\mu}(x)(y)= y-{{1}\over {\mu}}T'(x)(y)\ ,$$
the above recalled characterization implies that the functional $I_{\mu}$ is not convex.
Since $T$ is compact, on the one hand, $J_T$ is sequentially
weakly continuous ([29], Corollary 41.9) and, on the other hand, in view of $(7.2)$ the operator
$I_{\mu}'$ (recall that $-{{1}\over {\mu}}>\lambda_0$) is proper ([28], Example 4.43).
 The compactness of $T$ also implies that, for each
$x\in X$, the operator $T'(x)$ is compact ([28], Proposition 7.33) and so, for each $\lambda\in {\bf R}$,
 the operator $y\to y+\lambda T'(x)(y)$ is Fredholm of index zero ([28], Example 8.16).
Therefore, the operator
 $I_{\mu}'$ is non-monotone, proper and Fredholm of index zero.
Clearly, by $(7.1)$, the functional $x\to I_{\mu}(x)+\langle z,x\rangle$ is coercive for all $z\in X$. Then, in view of 
Corollary 5.1, the operator $I_{\mu}'$ is not injective.  At this point, we can apply Theorem 7.B
to infer that the set $\tilde S_{I_{\mu}'}$ contains an accumulation point. Finally, notice that
$$\tilde S_{I_{\mu}'}=A_{\mu}\ ,$$
and $(iii)$ follows.
The implication $(iii)\to (ii)$ is trivial. Finally, the implication $(ii)\to (i)$ is
provided by Theorem 2.1.11 of [27] again.\hfill $\bigtriangleup$

\bigskip
{\bf 8. Integral functionals on $L^p$-spaces}
\bigskip
In this section, we present an application of Theorem 3.1 to integral functionals on $L^p$-spaces.
The main general result is Theorem 8.1 below from which, in turn, we derive a series of consequences.\par
\medskip
In the sequel,  $(T,{\cal F},\mu)$ ($\mu(T)>0$) is a $\sigma$-finite
measure space, $Y$ is a reflexive real Banach space and
$\varphi, \psi:Y\to {\bf R}$ are two sequentially weakly lower semicontinuous
functionals such that
$$\inf_{y\in Y}{{\min\{\varphi(y),\psi(y)\}}\over
{1+\|y\|^p}}>-\infty \eqno{(8.1)}$$
for some $p>0$. \par
\smallskip
For each $\lambda\in [0,\infty]$, we denote by $M_{\lambda}$ the set of all
global minima of $\varphi+\lambda\psi$ or the empty set according to whether
$\lambda<+\infty$ or $\lambda=+\infty$. We adopt the conventions $\inf\emptyset=+\infty$
and $\sup\emptyset=-\infty$.\par
\smallskip
Moreover, $a, b$ are two fixed numbers in $[0,+\infty]$, with $a<b$, and
$\alpha$, $\beta$ are the
numbers so defined:
$$\alpha=\max\left\{\inf_Y\psi,\sup_{M_b}\psi\right\}\ ,$$
$$\beta=\min\left\{\sup_Y\psi,\inf_{M_a}\psi\right\}\ .$$
\smallskip
As usual, $L^p(T,Y)$ denotes
the space of all $\mu$-strongly measurable functions $u:T\to Y$ such that
$$\int_T\|u(t)\|^pd\mu<+\infty\ .$$
\medskip
THEOREM 8.1. - {\it  Assume that the functional
 $\varphi+\lambda\psi$
is coercive and has a unique global minimum for each $\lambda\in ]a,b[$.
 Assume also that
$$\alpha<\beta\ .$$\par
Then, for each
$\gamma\in L^{\infty}(T)\cap L^1(T)\setminus \{0\}$, with $\gamma\geq 0$,
and for each $r\in ]\alpha,\beta[$, if we put
$$V_{\gamma,r}=
\left \{u\in L^p(T,Y) : \int_T\gamma(t)\psi(u(t))d\mu\leq
r\int_T\gamma(t)d\mu\right
\}\ ,$$
we have
$$\inf_{u\in V_{\gamma,r}}
\int_T\gamma(t)\varphi(u(t))d\mu=
\inf_{\psi^{-1}(r)}\varphi\int_T\gamma(t)d\mu\ .\eqno{(8.2)}$$}\par
\smallskip
PROOF. First, we also assume that
$$\varphi(0)=\psi(0)=0\ .$$
Actually, once we prove the theorem under this additional assumption, 
the general version is obtained applying the particular version to the
functions $\varphi-\varphi(0)$ and $\psi-\psi(0)$.
Next, observe that the functionals
$\varphi$ and $\psi$ are Borel (in the weak topology, and so in the
strong one too). This implies that, for each $u\in L^p(T,Y)$,
the functions $\varphi\circ u$ and $\psi\circ u$ are $\mu$-measurable. On the
other hand,
in view of $(8.1)$, for some $c>0$, we have
$$-c\gamma(t)(1+\|u(t)\|^p)\leq \gamma(t)\min\{\varphi(u(t)),\psi(u(t))\}$$
for all $t\in T$. Since $\gamma\in L^{\infty}(T)\cap L^1(T)$,
the function $t\to -\gamma(t)(1+\|u(t)\|^p)$ lies in $L^1(T)$, and so
the integrals $\int_T\gamma(t)\varphi(u(t))d\mu$ and
$\int_T\gamma(t)\psi(u(t))d\mu$ exist and belong to $]-\infty,+\infty]$.
For each $\lambda\in ]a,b[$,
denote by $\hat y_{\lambda}$
the unique global minimum in $Y$ of the functional
$\varphi+\lambda\psi$.
 By Theorem 3.1 and Remark 3.1, 
there exists $\lambda_r\in ]a,b[$ such that
$$\psi(\hat y_{\lambda_r})=r\ .$$
So, we have
$$\varphi(\hat y_{\lambda_r})+\lambda_rr\leq \varphi(y)+\lambda_r\psi(y)$$
for all $y\in Y$. From this, it clearly follows that
$$\varphi(\hat y_{\lambda_r})=\inf_{\psi^{-1}(r)}\varphi\ . \eqno{(8.3)}$$
 Likewise,
for each $u\in L^p(T,Y)$, it follows that
$$(\varphi(\hat y_{\lambda_r})+\lambda_rr)\int_T\gamma(t)d\mu\leq
\int_T(\gamma(t)(\varphi(u(t))+\lambda_r\psi(u(t)))d\mu\ .$$
Therefore, for each $u\in V_{\gamma,r}$, we have
$$\varphi(\hat y_{\lambda_r})\int_T\gamma(t)d\mu\leq
\int_T\gamma(t)\varphi(u(t))d\mu\ ,$$
and hence
$$\varphi(\hat y_{\lambda_r})\int_T\gamma(t)d\mu\leq
\inf_{u\in V_{\gamma,r}}\int_T\gamma(t)\varphi(u(t))d\mu\ .\eqno{(8.4)}$$
In view of $(8.3)$, to get $(8.2)$, we have to show that
$$\varphi(\hat y_{\lambda_r})\int_T\gamma(t)d\mu=
\inf_{u\in V_{\gamma,r}}\int_T\gamma(t)\varphi(u(t))d\mu\ .\eqno{(8.5)}$$
Arguing by contradiction, assume that $(8.5)$ does not hold. So, in
view of $(8.4)$, we would have
$$\varphi(\hat y_{\lambda_r})\int_T\gamma(t)d\mu<
\inf_{u\in V_{\gamma,r}}\int_T\gamma(t)\varphi(u(t))d\mu\ .\eqno{(8.6)}$$
 From $(8.6)$, in turn, as
 $(T,{\cal F},\mu)$ is $\sigma$-finite, it would follow the existence of
$\tilde T\in {\cal F}$, with $\mu(\tilde T)<+\infty$, such that
$$\varphi(\hat y_{\lambda_r})\int_{\tilde T}\gamma(t)d\mu<
\inf_{u\in V_{\gamma,r}}\int_T\gamma(t)\varphi(u(t))d\mu\ .\eqno{(8.7)}$$
Now, consider the function $\hat u:T\to Y$ defined by
$$\hat u(t)=\cases {\hat y_{\lambda_r} & if $x\in \tilde T$\cr & \cr
0 & if $x\in T\setminus\tilde T$\ .\cr}$$
Clearly, $\hat u\in L^p(T,Y)$. We also have
$$\int_T\gamma(t)\psi(\hat u(t))d\mu=
\int_{\tilde T}\gamma(t)\psi(\hat u(t))d\mu\leq r\int_T\gamma(t)d\mu$$
and so $\hat u\in V_{\gamma,r}$. But
$$\int_T\gamma(t)\varphi(\hat u(t))d\mu=\varphi(\hat y_{\lambda_r})
\int_{\tilde T}\gamma(t)d\mu$$
and this contradicts $(8.7)$. The proof is complete.\hfill $\bigtriangleup$\par
\medskip
REMARK 8.1. -  In general, the conclusion of Theorem 8.1 is no longer true if,
for some $\lambda\in ]a,b[$,
the function $\varphi+\lambda\psi$
has more than one global minimum. A simple example (with $a=0$ and $b=+\infty$)
is provided by taking
$Y={\bf R}$,
$$\varphi(y)=\cases {y^2 & if $y\leq 1$\cr & \cr
2-y & if $y>1$\cr}$$
and
$$\psi(y)=y^2\ .$$
So, $\varphi$ is unbounded below and $\varphi+\lambda\psi$ is coercive for all $\lambda>0$.
Clearly, we have $\alpha=0$ and $\beta=+\infty$.
However, for $r=1$, $(8.2)$ is not satisfied, since $0\in V_{\gamma,r}$,
$\int_T\gamma(t)\varphi(0)d\mu=0$, while $\inf_{\psi^{-1}(1)}\varphi=1$.\par
\medskip
REMARK 8.2. - At present, we do not know if the conclusion of Theorem 8.1 holds 
without the coercivity assumption on $\varphi+\lambda\psi$.
\par
\medskip
We now consider a series of consequences of Theorem 8.1.\par
\smallskip
First, we want to state explicitly the form that Theorem 8.1 assumes when
$T={\bf N}$, ${\cal F}$ is the power set of ${\bf N}$ and
$$\mu(A)=\hbox {\rm card}(A)$$
for all $A\subseteq {\bf N}$.\par
\smallskip
  Denote by $l_p(Y)$ the
space of all sequences $\{u_n\}$ in $Y$ such that
$$\sum_{n=1}^{\infty}\|u_n\|^p<+\infty\ .$$
\medskip
THEOREM 8.2. - {\it Let $\varphi, \psi$ satisfy the assumptions of Theorem 8.1.
\par
Then, for each sequence $\{a_n\}\in l_1({\bf R})\setminus
\{0\}$, with $\inf_{n\in {\bf N}}
a_n\geq 0$,
and for each $r\in ]\alpha,\beta[$, if we put
$$V_{\{a_n\},r}=
\left \{\{u_n\}\in l_p(Y) : \sum_{n=1}^{\infty}a_n\psi(u_n)\leq
r\sum_{n=1}^{\infty}a_n\right\}\ ,$$
we have
$$\inf_{\{u_n\}\in V_{\{a_n\},r}}
\sum_{n=1}^{\infty}a_n\varphi(u_n)=
\inf_{\psi^{-1}(r)}\varphi\sum_{n=1}^{\infty}a_n\ .$$}
\medskip
The next two results deals with consequences of Theorem 8.1 in the case where
$\varphi\in Y^{*}\setminus \{0\}$.\par
\medskip
THEOREM 8.3. - {\it Let $y\to \|y\|^q$ be strictly convex for some $q>1$
 and let $\varphi$ be non-zero, continuous and linear.
Moreover, let $\eta:[0,+\infty[\to {\bf R}$ be an increasing strictly
convex function.\par
Then, for each
$\gamma\in L^{\infty}(T)\cap L^1(T)\setminus \{0\}$, with $\gamma\geq 0$,
and for each $r>\eta(0)$ and $p\geq 1$, if we put
$$V_{\gamma,r}=
\left \{u\in L^p(T,Y) : \int_T\gamma(t)\eta(\|u(t)\|^q)d\mu\leq
r\int_T\gamma(t)d\mu\right
\}\ ,$$
we have
$$\inf_{u\in V_{\gamma,r}}
\int_T\gamma(t)\varphi(u(t))d\mu=
-\|\varphi\|_{Y^*}(\eta^{-1}(r))^{1\over q}\int_T\gamma(t)d\mu\ .$$}
\smallskip
PROOF.
 By the assumptions made on $\eta$, the functional
 $y\to \eta(\|y\|^q)$ is strictly convex and, for some $m, c>0$, one
 has
 $$\eta(t)\geq mt-c$$
 for all $t\geq 0$. As a consequence, for each $\lambda>0$, the functional
$y\to \varphi(y)+\lambda\eta(\|y\|^q)$ is coercive and has a unique
global minimum in $X$. At this point, the conclusion follows directly
from Theorem 8.1, applied taking $a=0$, $b=+\infty$,
$\psi(y)=\eta(\|y\|^q)$ and observing that
$(8.1)$ holds for each $p\geq 1$ and that $\alpha=\eta(0)$ , $\beta=+\infty$.
\hfill $\bigtriangleup$\par
\medskip
THEOREM 8.4. - {\it Let $\varphi$ be non-zero, continuous and linear and
let $\psi$ be $C^1$ with
$$\lim_{\|y\|\to +\infty}{{\psi(y)}\over {\|y\|}}=+\infty\ .\eqno{(8.8)}$$
Finally, assume that, for each $\mu<0$, the equation
$$\psi'(y)=\mu\varphi \eqno{(8.9)}$$ 
has a unique solution in $Y$ or even at most two when $\hbox {\rm dim}(Y)<
\infty$\ .\par
Then, for each $p\geq 1$, the conclusion of Theorem 8.1 holds with
any $r>\inf_Y\psi$\ .}\par
\smallskip
PROOF.  In view of $(8.8)$, the functional $\varphi+\lambda
\psi$ is coercive for each $\lambda>0$. Let $\hat x$ be a global minimum
of this functional. Then, $\hat x$ satisfies $(8.9)$ with $\mu=-\lambda^{-1}$.
So, when dim$(Y)=\infty$, the uniqueness of $\hat x$ follows from an
assumption directly. Now, assume that dim$(Y)<\infty$. In this case,
$\varphi+\lambda\psi$ satisfies the Palais-Smale condition. As a consequence,
if $\varphi+\lambda\psi$ was admitting two global minima, then,
thanks to Corollary 1 of [9], $(8.9)$ would have at least three solutions
for $\mu=-\lambda^{-1}$, against an assumption. Now, we can apply
Theorem 8.1, with $p\geq 1$, $a=0$, $b=+\infty$, observing that
$\alpha=\inf_Y\psi$ and $\beta=+\infty$. \hfill
$\bigtriangleup$\par
\medskip
Here is a consequence of Theorem 8.1 in the case when $Y$ is a Hilbert
space and $\varphi$ has a Lipschitzian derivative:\par
\medskip
THEOREM 8.5. - {\it Let $Y$ be a Hilbert space,
let $\varphi$ be $C^1$ and let $\varphi'$ be Lipschitzian,
with Lipschitz constant $L>0$. Assume that $\varphi'(0)\neq 0$.
Set
$$S=\{y\in Y: \varphi'(y)+L y=0\}$$
and
$$\rho=\inf_{y\in S}\|y\|^2\ .$$
Then, for each
$\gamma\in L^{\infty}(T)\cap L^1(T)\setminus \{0\}$,
with $\gamma\geq 0$,
and for each $r\in ]0,\rho[$, $p\geq 2$, if we put
$$V_{\gamma,r}=
\left \{u\in L^p(T,Y) : \int_T\gamma(t)\|u(t)\|^2d\mu\leq
r\int_T\gamma(t)d\mu\right
\}\ ,$$
we have
$$\inf_{u\in V_{\gamma,r}}\int_T\gamma(t)
\varphi(u(t))d\mu=\inf_{\|y\|^2=r}\varphi(y)\int_T\gamma(t)d\mu
\ .$$}
\smallskip
PROOF. 
Note that the functional $y\to \varphi(y)+{{\lambda}\over {2}}\|y\|^2$
is convex if $\lambda=L$, while it is strictly convex and coercive
if $\lambda>L$  (see, for instance, Proposition 2.2 of [16]). So, this functional
has a unique global minimum if
$\lambda>L$, while the set of its global minima coincides with $S$ if
$\lambda=L$.
At this point, the conclusion is obtained applying Theorem 8.1
with
$$\psi(y)={{\|y\|^2}\over {2}}$$
for all $y\in Y$ and
$$a=L\ , b=+\infty\ ,$$
taking into account that $(8.1)$ is satisfied for each $p\geq 2$ since
$\varphi'$ is Lipschitzian and observing that $\alpha=0$ and $\beta={{\rho}\over
{2}}$.
\hfill $\bigtriangleup$
\medskip
In the next result, we will apply Theorem 8.1 taking as $Y$ the usual
Sobolev space $W^{1,q}_0(\Omega)$ with the usual norm
$$\left ( \int_{\Omega}|\nabla v(x)|^qdx\right )^{1\over q} ,$$
where $\Omega$ is bounded domain in ${\bf R}^n$ ($n\geq 3$)
 with smooth boundary and $q>1$.\par
\smallskip
Moreover, if $u\in L^p(T,W^{1,q}_0(\Omega))$ we will write $u(t,x)$ instead
of $u(t)(x)$. That is, we will identify $u$ with
the function $\omega:T\times \Omega\to {\bf R}$ defined by
$$\omega(t,x)=u(t)(x)$$
for all $(t,x)\in T\times\Omega$.  \par
\medskip
THEOREM 8.6. - {\it Let $f:{\bf R}\to [0,+\infty[$ be a continuous
function, with $f(0)=0$ and
$\liminf_{\xi\to +\infty}f(\xi)>0$,
 such that $\xi\to {{f(\xi)}\over {\xi^{q-1}}}$
is decreasing in $]0,+\infty[$ and
$$\lim_{|\xi|\to +\infty}{{f(\xi)}\over {|\xi|^{q-1}}}=0 \eqno{(8.10)}$$
for some $q>1$.\par
Then, for each
$\gamma\in L^{\infty}(T)\cap L^1(T)\setminus \{0\}$, with $\gamma\geq 0$,
and each $r>0$, $p\geq q$, if we put
$$V_{\gamma,r}=\left \{u\in L^p(T,W^{1,q}_0(\Omega)) : \int_T\gamma(t)\left (
\int_{\Omega}|\nabla u(t,x)|^qdx\right ) d\mu\leq r\int_T\gamma(t)d\mu\right \}\ ,$$
we have
$$\sup_{u\in V_{\gamma,r}}\int_T\gamma(t)\left ( \int_{\Omega}F(u(t,x))dx\right )d\mu=
\sup_{v\in W^{1,q}_0(\Omega), \int_{\Omega}|\nabla v(x)|^qdx=r}
\int_{\Omega}F(v(x))dx\int_T\gamma(t)d\mu\ ,$$
where 
$$F(\xi)=\int_0^{\xi}f(s)ds$$
for all $\xi\in {\bf R}$.} \par
\smallskip
PROOF. We are going to apply Theorem 8.1 taking $Y=W^{1,q}_0(\Omega)$ and
$$\varphi(v)=-\int_{\Omega}F(v(x))dx\ ,$$
$$\psi(v)=\int_{\Omega}|\nabla v(x)|^qdx$$
for all $v\in W^{1,q}_0(\Omega)$. Due to $(8.10)$,
 by classical results, 
$\varphi$ is sequentially weakly continuous in $W^{1,q}_0(\Omega)$,
$(8.1)$ is satisfied for any $p\geq q$, and, for each $\lambda>0$,
the functional $\varphi+\lambda\psi$ is $C^1$, coercive and
satisfies the Palais-Smale condition. Moreover, since $f\geq 0$,
its non-zero critical points are strictly positive in $\Omega$ ([3], [26]).
 Moreover, since
the function $\xi\to {{f(\xi)}\over {\xi^{q-1}}}$
is decreasing in $]0,+\infty[$, Proposition 4.2 of [4] ensures that,
for each $\lambda>0$, there exists at most one strictly
positive critical point of $\varphi+\lambda\psi$. As a consequence,
we infer that, for each $\lambda>0$, the functional $\varphi+\lambda\psi$
has a unique global minimum in $W^{1,q}_0(\Omega)$, since otherwise,
in view of Corollary 1 of [9], it would have at least three critical
points. Hence, we are allowed to apply Theorem 8.1 with $a=0$ and
$b=+\infty$. Clearly, we have $\alpha=0$ and $\beta=+\infty$ (since
$\lim_{\xi\to +\infty}F(\xi)=+\infty$ and hence $\varphi$ is
unbounded below). The proof is complete.\hfill $\bigtriangleup$\par
\medskip
The next application of Theorem 8.1 concerns
a Jensen-like inequality in $L^p$-spaces.\par
\medskip
THEOREM 8.7. - {\it Let $f:{\bf R}\to {\bf R}$
 be a continuous function, positive and differentiable in $]0,+\infty[$,
with $\sup_{]-\infty,0]}f\leq 0$. Assume that, for some
$\delta\geq 0$, the
function $y\to \delta|y|^p-f(y)$ has no global minima in
${\bf R}$, 
$$\limsup_{y\to +\infty}{{f(y)}\over {y^p}}=\delta \eqno{(8.11)}$$
and the function
$$y\to {{f'(y)}\over {y^{p-1}}}$$
 is injective in $]0,+\infty[$.\par
Then, for each $\gamma\in L^{\infty}(T)\cap L^1(T)\setminus \{0\}$,
with $\gamma\geq 0$, 
one has
$$\int_T\gamma(t)f(u(t))d\mu\leq
f\left (\left ( {{\int_T\gamma(t)|u(t)|^pd\mu}\over {\int_T\gamma(t)d\mu}}\right )
^{1\over p}\right )\int_T\gamma(t)d\mu\ ,$$
for all $u\in L^p(T)$\ .}
\smallskip
PROOF. We are going to apply Theorem 8.1 with $Y={\bf R}$, $\varphi(y)=
-f(y)$, $\psi(y)=|y|^p$ and $a=\delta$, $b=+\infty$. 
 Fix
$\lambda>\delta$.
 From $(8.11)$, we clearly infer that $\varphi+\lambda\psi$ is coercive. We
now show that this function has a unique global minimum. Arguing by
contradiction, assume that $y_1, y_2\in {\bf R}$ are two distinct
global minima of $\varphi+\lambda\psi$. We can suppose that $y_1<y_2$.
Since $\varphi(y)+\lambda\psi(y)>0$ for all $y<0$ and $\varphi(0)+\lambda\psi(0)=0$,
it would follow that $y_1\geq 0$. By the Rolle theorem, there would be
$y_3\in ]y_1,y_2[$ such that 
$$p\lambda y_3^{p-1}=f'(y_3)\ .$$
As a consequence, we would have
$${{f'(y_2)}\over {y_2^{p-1}}}={{f'(y_3)}\over {y_3^{p-1}}}\ ,$$
contrary to the assumption that the function $y\to {{f'(y)}\over {y^{p-1}}}$ is
injective in $]0,+\infty[$.
So, we are allowed to apply Theorem 8.1, observing that
$\alpha=0$ and $\beta=+\infty$. Let $u\in L^p(T)$. Put
$$r={{\int_T\gamma(t)|u(t)|^pd\mu}\over {\int_T\gamma(t)d\mu}}\ .$$
If $r=0$ the inequality to prove is clear: both sides are zero. So, assume $r>0$.
Clearly, we have
$$\inf_{\psi^{-1}(r)}\varphi=-f\left (\left ( {{\int_T\gamma(t)|u(t)|^pd\mu}\over {\int_T\gamma(t)d\mu}}\right )
^{1\over p}\right )$$
and hence, since $u\in V_{\gamma,r}$, it follows
$$\int_T\gamma(t)f(u(t))d\mu\leq
f\left (\left ( {{\int_T\gamma(t)|u(t)|^pd\mu}\over {\int_T\gamma(t)d\mu}}\right )
^{1\over p}\right )\int_T\gamma(t)d\mu\ ,$$
as claimed.\hfill$\bigtriangleup$
\medskip
REMARK 8.3. -
The class of functions $f$ satisfying the assumptions of Theorem 8.7 is quite
broad. For instance, a typical function in that class is
$$f(y)=a_0\log(1+(y^+)^p)+\sum_{i=1}^k a_i (y^+)^{q_i}$$
where $y^+=\max\{y,0\}$,
$a_i$ ($i=0,...,k$) are $k+1$ non-negative numbers, with
$\sum_{i=0}^k a_i>0$, and $q_i$  ($i=1,...,k$)
are $k$ positive numbers less than $p$.
\medskip
As a consequence of this remark, we get, for instance, the following\par
\medskip
COROLLARY 8.1. - {\it For each $\gamma\in L^{\infty}(T)\cap L^1(T)\setminus
\{0\}$, with $\gamma\geq 0$, one has
$$\int_T\gamma(t)\log(1+|u(t)|^p)d\mu\leq
\log\left ( 1+{{\int_T\gamma(t)|u(t)|^pd\mu}\over {\int_T\gamma(t)
d\mu}}
\right )\int_T\gamma(t)d\mu \eqno{(8.12)}$$
for all $u\in L^p(T)$.}\par
\bigskip
{\bf 9. Integral functionals on Sobolev spaces}\par
\bigskip
From now on, $\Omega\subset {\bf R}^n$ is a bounded domain with smooth boundary.\par
\smallskip
For $p>1$, on the Sobolev space $W^{1,p}(\Omega)$ we consider the norm
$$\|u\|=\left ( \int_{\Omega}|\nabla
u(x)|^{p}dx+
\int_{\Omega}|u(x)|^{p}dx\right ) ^{{1}\over {p}}\ .$$
If $p>n$, $W^{1,p}(\Omega)$
is compactly embedded in $C^0(\overline {\Omega})$ and hence
the constant
$$c=\sup_{u\in W^{1,p}(\Omega)\setminus \{0\}}
{{\sup_{x\in \Omega}|u(x)|}\over {\|u\|}} \eqno{(9.1)}$$
is finite.\par
\smallskip
Recall that a function $f:\Omega\times {\bf R}^{m}\to ]-\infty,+\infty]$
 is said to be a normal integrand ([22]) if it is
${\cal L}(\Omega)\otimes {\cal B}({\bf R}^m)$-measurable and
$f(x,\cdot)$ is lower semicontinuous for a.e. $x\in \Omega$.
Here ${\cal L}(\Omega)$ and ${\cal B}({\bf R}^m)$ denote the
Lebesgue and the Borel $\sigma$-algebras of subsets of $\Omega$ and
${\bf R}^m$, respectively.  \par
\smallskip
Recall that if $f$ is a normal integrand, then, for each measurable function $u:\Omega\to
{\bf R}^m$, the composite function $x\to f(x,u(x))$ is measurable
([22]).
\par
\smallskip
If $\xi\in {\bf R}$, we continue to denote by $\xi$ the constant function
on $\Omega$ assuming the value $\xi$.\par
\smallskip 
The next result is an application of Theorem 5.1.\par
\medskip
THEOREM 9.1. - {\it Assume $p>n$. Let $f:\Omega\times {\bf R}\to ]-\infty,+\infty]$
and $\varphi:\Omega\times {\bf R}\times {\bf R}^n\to ]-\infty,+\infty]$ be
two normal integrands satisfying the following conditions:\par
\noindent
$(i)$\hskip 5pt there exist $\nu>0$ and $\gamma\in L^1(\Omega)$ such that
$$\nu(|\xi|^p+|\eta|^p)+\gamma(x)\leq \varphi(x,\xi,\eta)$$
for all $(x,\xi,\eta)\in \Omega\times {\bf R}\times {\bf R}^n$, and,
for each $(x,\xi)\in \Omega\times {\bf R}$, the function
$\varphi(x,\xi,\cdot)$ is convex in ${\bf R}^n$\ ;\par
\noindent
$(ii)$\hskip 5pt for each $\epsilon>0$, there exists $\gamma_{\epsilon}\in L^1(\Omega)$ such that
$$-\epsilon|\xi|^p+\gamma_{\epsilon}(x)\leq f(x,\xi)$$
for all $(x,\xi)\in \Omega\times {\bf R}$\ ;
\par
\noindent
$(iii)$\hskip 5pt there exist $\xi_1, \rho\in {\bf R}$ such that
$$\int_{\Omega}\varphi(x,\xi_1,0)dx<\rho\ ,\hskip 5pt
\int_{\Omega}f(x,\xi_1)dx<+\infty$$
and
$$f(x,\xi_1)=\inf_{|\xi|\leq \delta}f(x,\xi)$$
for all $x\in \Omega$,
where
$$\delta=
c\left ( {{\rho-\int_{\Omega}\gamma(x)dx}\over
{\nu}}\right ) ^{1\over p}$$
and $c$ is given in $(9.1)$.\par
Under such hypotheses, for every sequentially weakly closed
set $V\subseteq W^{1,p}(\Omega)$ containing the constant $\xi_1$
and a $w$ for which
$$\int_{\Omega}\varphi(x,w(x),\nabla w(x))dx<+\infty$$
and
$$\int_{\Omega}f(x,w(x))dx<\int_{\Omega}f(x,\xi_1)dx\ ,$$
there exists $\lambda^*>0$ such that the restriction to $V$ of the
functional 
$$u\to \int_{\Omega}f(x,u(x))dx+\lambda^*\int_{\Omega}
\varphi(x,u(x),\nabla u(x))dx$$  has at least two global
 minima.}\par
\smallskip
PROOF. For each $u\in W^{1,p}(\Omega)$, set
$$\tilde J(u)=\int_{\Omega}f(x,u(x))dx$$
and
$$\tilde\Phi(u)=\int_{\Omega}\varphi(x,u(x),\nabla u(x))dx\ .$$
By a classical result ([2], Theorem 4.6.8), for each $\lambda>0$
the functional $\tilde J+\lambda \tilde\Phi$ is sequentially weakly
lower semicontinuous. On the other hand, for $\epsilon\in
]0,\lambda\nu[$, by $(ii)$, we have
$$\tilde J(u)+\lambda\tilde\Phi(u)\geq (\lambda\nu-\epsilon)\|u\|+\int_{\Omega}\gamma_{\epsilon}(x)dx\ .$$
Consequently, by reflexivity and Eberlein-Smulyan theorem, the
sub-level sets of $\tilde J+\lambda\tilde\Phi$ are weakly compact. 
Now, let $V\subseteq W^{1,p}(\Omega)$ be as in the conclusion.
Set
$$X=\{u\in V : \sup\{\tilde J(u),\tilde\Phi(u)\}<+\infty\}\ .$$
Observe that $\xi_1, w\in X$ and that
$$\{u\in X:\tilde J(u)+\lambda\tilde\Phi(u)\leq r\}=
\{u\in V:\tilde J(u)+\lambda\tilde\Phi(u)\leq r\} \eqno{(9.2)}$$
for all $\lambda>0$, $r\in {\bf R}$.
Denote by $J$ and $\Phi$ the restrictions to $X$ of $\tilde J$ and $\tilde\Phi$
respectively. We want to apply Theorem 5.1 considering $X$ with the relative
weak topology. Clearly, in view of $(9.2)$, $(a_1)$ holds.
Concerning $(b_1)$,
observe that for each $u\in \Phi^{-1}(]-\infty,\rho])$, by $(i)$, one has
$$\nu\|u\|^p+\int_{\Omega}\gamma(x)dx\leq \rho$$
and so
$$\sup_{\Omega}|u|\leq c\left ( {{\rho-\int_{\Omega}\gamma(x)dx}\over
{\nu}}\right ) ^{1\over p}\ ,$$
the above inequalities being strict if $\Phi(u)<\rho$. Then, from this and from $(iii)$, it follows that
$$J(\xi_1)=\inf_{\Phi^{-1}(]-\infty,\rho])}J$$
and
$$\Phi(\xi_1)<\rho$$
as well.
Consequently, $(b_1)$ is satisfied taking $u_1=\xi_1$ and
$u_2=w$.
So, the conclusion follows directly from Theorem 5.1.
\hfill $\bigtriangleup$\par
\medskip
Let $p>1$.
If $n\geq p$, we denote by ${\cal A}_p$ the class of all
continuous functions $f:{\bf R}\to {\bf R}$ such that
$$\sup_{\xi\in {\bf R}}{{|f(\xi)|}\over
{1+|\xi|^{s}}}<+\infty\ ,$$
where  $0<s< {{pn-n+p}\over {n-p}}$ if $p<n$ and $0<s<+\infty$ if
$p=n$. While, when $n<p$, ${\cal A}$ stands for the class
of all continuous functions $f:{\bf R}\to {\bf R}$. 
Given $f\in {\cal A}_p$, consider the following Dirichlet problem
$$\cases {-\hbox {\rm div}(|\nabla u|^{p-2}\nabla u)=
f(u)
 & in
$\Omega$\cr & \cr u=0 & on
$\partial \Omega$\ .\cr}\eqno{(P_{f})} $$
 Let us recall
that a weak solution
of $(P_{f})$ is any $u\in W^{1,p}_0(\Omega)$ such that
 $$\int_{\Omega}|\nabla u(x)|^{p-2}\nabla u(x)\nabla v(x)dx
-\int_{\Omega}
f(u(x))v(x)dx=0$$
for all $v\in W^{1,p}_0(\Omega)$.\par
\smallskip
Moreover, $\lambda_{1,p}$ denotes the principal
eigenvalue of the problem
$$\cases {-\hbox {\rm div}(|\nabla u|^{p-2}\nabla u)=
\lambda  |u|^{p-2}u
 & in
$\Omega$\cr & \cr u=0 & on
$\partial \Omega$\ .\cr} $$
We have
$$\lambda_{1,p}=\inf_{u\in W^{1,p}_0(\Omega)\setminus \{0\}}{{\int_{\Omega}|\nabla u(x)|^pdx}\over
{\int_{\Omega}|u(x)|^pdx}}\ .$$
\medskip
Also, let us recall the following consequence of the variational principle established in [10]:\par
\medskip
THEOREM 9.A. - {\it Let $X$ be a reflexive real Banach space and let
$\Phi, \Psi:X\to {\bf R}$ be two sequentially weakly lower semicontinuous
functionals, with $\Phi(0)=\Psi(0)=0$, and with $\Psi$ also coercive and continuous.\par
Then, for each $\sigma>\inf_X\Psi$ and each $\lambda$ satisfying
$$\lambda>-{{\inf_{\Psi^{-1}(]-\infty,\sigma])}\Phi}\over {\sigma}}$$
the functional $\lambda\Psi+\Phi$  has a local
minimum belonging to $\Psi^{-1}(]-\infty,\sigma[)$\ .}\par
\medskip
The next result is an application of Theorem 8.1.\par
\medskip
THEOREM 9.2. - {\it Let $f\in {\cal A}_p$, with $f\geq 0$, 
 and let $F(\xi)=\int_0^{\xi}F(t)dt$ for all $\xi\in {\bf R}$. 
Assume that:\par
\noindent
$(a_1)$\hskip 5pt $\lim_{\xi\to 0^+}{{F(\xi)}\over {\xi^p}}=+\infty\ ;$\par
\noindent
$(a_2)$\hskip 5pt $\delta:=\limsup_{\xi\to +\infty}{{F(\xi)}\over {\xi^p}}<+\infty\ ;$\par
\noindent
$(a_3)$\hskip 5pt the function $\xi\to \delta \xi^p-F(\xi)$ has no global minima in $[0,+\infty[$\ ;\par
\noindent
$(a_4)$\hskip 5pt 
for each $\lambda>p\delta$, the equation $\lambda \xi^{p-1}=f(\xi)$ has
at most two solutions in $]0,+\infty[$\ .\par
\noindent
Under such hypotheses, for each $\rho>0$ and each $\nu\in ]0,1]$ satisfying
$$\nu<{{\lambda_{1,p}\rho^p}\over {p F(\rho)}}\ ,\eqno{(9.3)}$$
  the problem
$$\cases {-\hbox {\rm div}(|\nabla u|^{p-2}\nabla u)=
\nu f(u)
 & in
$\Omega$\cr & \cr u=0 & on
$\partial \Omega$\cr} $$
has a positive weak solution satisfying
$$\int_{\Omega}|\nabla u(x)|^pdx<\rho^p\lambda_{1,p}\hbox {\rm meas}(\Omega)\ .$$}
\smallskip
PROOF. Fix $\rho$ and $\nu$ as above.  
Since $f\geq 0$, by classical results ([3], [27]), the positive weak solutions of
the problem are exactly the non-zero critical points in $W^{1,p}_0(\Omega)$ of the
energy functional 
$$u\to {{1}\over {p}}\int_{\Omega}|\nabla u(x)|^pdx-\nu\int_{\Omega}F(u(x))dx\ .$$
 We are going to apply Theorem 8.1 taking $Y={\bf R}$, $\varphi(\xi)=-\nu F(\xi)$, $\psi(\xi)=|\xi|^p$,
$a=\delta$ and $b=+\infty$. Note that $\varphi$ is
non-negative in $]-\infty,0]$. So, $(8.1)$ is satisfied in view of $(a_2)$.
Fix $\lambda>\delta$.
From $(a_2)$ again, it follows that $\varphi+\lambda\psi$ is coercive . Arguing by contradiction,
assume that $\varphi+\lambda\psi$ has two global minima, say $\xi_1, \xi_2$, with
$\xi_1<\xi_2$. From $(a_1)$ it
follows that 
$$\inf_{[0,+\infty[}(\varphi+\lambda\psi)<0\ .$$
This fact implies that $\xi_1>0$. As a consequence,
 the equation
$$p\lambda \xi^{p-1}=\nu f(\xi)$$
would admit the solutions $\xi_1$, $\xi_2$ and a third one in $]\xi_1,\xi_2[$ given by
the Rolle theorem. But, this contradicts $(a_4)$ .
 Hence, the function
$\varphi+\lambda\psi$ has a unique global minimum. Further, note that $\alpha=0$ 
and, in view of $(a_3)$, $\beta=+\infty$. Then, if we put
$$V_{\rho}=\left\{u\in L^p(\Omega) : \int_{\Omega}|u(x)|^pdx\leq \rho^p\hbox {\rm meas}(\Omega)\right\}\ ,$$
Theorem 8.1 ensures that
$$\sup_{u\in V_{\rho}}\int_{\Omega}F(u(x))dx=F(\rho)\hbox {\rm meas}(\Omega)\ .\eqno{(9.4)}$$
On the other hand, setting
$$B_{\rho}=\left\{u\in W^{1,p}_0(\Omega) : \int_{\Omega}|\nabla u(x)|^pdx\leq 
\rho^p\lambda_{1,p}\hbox {\rm meas}(\Omega)\right\}\ ,$$
we have
$$B_{\rho}\subseteq V_{\rho}\ .$$
Consequently
$$\sup_{u\in B_{\rho}}\int_{\Omega}F(u(x))dx\leq \sup_{u\in V_{\rho}}\int_{\Omega}F(u(x))dx\ .\eqno{(9.5)}$$
Now, if we put
$$\sigma=\rho^p\lambda_{1,p}\hbox {\rm meas}(\Omega)\ ,$$
in view of $(9.3)$ , $(9.4)$ and $(9.5)$, we have
$$\sup_{u\in W^{1,p}_0(\Omega), \int_{\Omega}|\nabla u(x)|^pdx\leq \sigma}
\int_{\Omega}\nu F(u(x))dx<{{\sigma}\over {p}}\ .$$
At this point, we can apply Theorem 9.A taking $X=W^{1,p}_0(\Omega)$,
$\Psi(u)=\int_{\Omega}|\nabla u(x)|^pdx$ and $\Phi(u)=-\nu\int_{\Omega}F(u(x))dx$.
Hence, the energy functional has a local minimum  $u$ (which is therefore a solution of the problem)  such that 
$$\int_{\Omega}|\nabla u(x)|^pdx<\rho^p\lambda_{1,p}\hbox {\rm meas}(\Omega)\ .$$
To show that $u\neq 0$, we finally remark that $0$ is not a local minimum of the
energy functional. Indeed, by a classical result, there is a
bounded and positive function $v\in W^{1,p}_0(\Omega)$ such that
$$\int_{\Omega}|\nabla v(x)|^pdx=\lambda_{1,p}\int_{\Omega}|v(x)|^pdx\ .$$
By $(a_1)$, there is $\theta>0$ such that
$$F(\xi)>{{\lambda_{1,p}}\over {\nu p}}\xi^p$$
for all $\xi\in ]0,\theta[$. Hence, for each $\eta\in \left ] 0, {{\theta}\over {\sup_{\Omega}v}}\right[$,
we have
$$\nu\int_{\Omega}F(\eta v(x))dx>{{\lambda_{1,p}}\over {p}}\int_{\Omega}|\eta v(x)|^pdx=
{{1}\over {p}}\int_{\Omega}|\nabla \eta v(x)|^pdx\ .$$ 
This shows that the energy functional takes negative values in each ball of $W^{1,p}_0(\Omega)$
centered at $0$ and so $0$ is not a local minimum for it. The proof is complete. \hfill
$\bigtriangleup$\par
\medskip
Note the following corollary of Theorem 9.2 (for the uniqueness, consider again Proposition 4.2 of [4]):\par
\medskip
COROLLARY 9.1. - {\it For each $\nu\in ]0,1]$, the unique positive weak solution of the problem
$$\cases {-\hbox {\rm div}(|\nabla u|\nabla u)=
\nu u
 & in
$\Omega$\cr & \cr u=0 & on
$\partial \Omega$\cr} $$
satisfies the inequality
$$\int_{\Omega}|\nabla u(x)|^3dx\leq {{27\hbox {\rm meas}(\Omega)}\over {8\lambda_{1,3}^2}}\nu^3\ .$$}
\medskip
Now, let $a, b\in {\bf R}$, with $a\geq 0$ and $b>0$.\par
\par
Consider the non-local problem
$$\cases {-\left ( a+b\int_{\Omega}|\nabla u(x)|^2dx\right )\Delta u =h(x,u) & in $\Omega$\cr
& \cr u=0 & on $\partial\Omega$\ ,\cr}$$
$h:\Omega\times {\bf R}\to {\bf R}$ being a Carath\'eodory function.\par
\smallskip
On the Sobolev space $H^1_0(\Omega)$, we consider the scalar product
$$\langle u,v\rangle=\int_{\Omega}\nabla u(x)\nabla v(x)dx$$
and the induced norm
$$\|u\|=\left ( \int_{\Omega}|\nabla u(x)|^2dx\right )^{1\over 2}\ .$$
We denote by ${\cal A}$ the class of all Carath\'eodory functions $f:\Omega\times {\bf R}\to {\bf R}$ such that
$$\sup_{(x,\xi)\in\Omega\times {\bf R}}{{|f(x,\xi)|}\over {1+|\xi|^p}}<+\infty \eqno{(9.6)}$$
for some $p\in \left ]0, {{n+2}\over {n-2}}\right [$.\par
\smallskip
Moreover, we denote by $\tilde{\cal A}$ the class of all Carath\'eodory functions $g:\Omega\times {\bf R}\to {\bf R}$ such that
$$\sup_{(x,\xi)\in\Omega\times {\bf R}}{{|g(x,\xi)|}\over {1+|\xi|^q}}<+\infty \eqno{(9.7)}$$
for some $q\in \left ] 0, {{2}\over {n-2}}\right [$.
\smallskip
Furthermore, we denote by $\hat {\cal A}$ the class of all functions $h:\Omega\times {\bf R}\to {\bf R}$ of the type
$$h(x,\xi)=f(x,\xi)+\alpha(x)g(x,\xi)$$
with $f\in {\cal A}, g\in\tilde{\cal A}$ and $\alpha\in L^2(\Omega)$.
For each $h\in\hat{\cal A}$, we define the functional $I_h:H^1_0(\Omega)\to {\bf R}$, by putting
$$I_f(u)=\int_{\Omega}H(x,u(x))dx$$
for all $u\in H^1_0(\Omega)$, where
$$H(x,\xi)=\int_0^{\xi}h(x,t)dt$$
for all $(x,\xi)\in \Omega\times {\bf R}$.\par
\smallskip
By classical results (involving the Sobolev embedding theorem), the functional $I_h$
turns out to be sequentially weakly continuous, of class $C^1$, with compact derivative given by
$$I_h'(u)(w)=\int_{\Omega}h(x,u(x))w(x)dx$$
for all $u,w\in H^1_0(\Omega)$.\par
\smallskip
Now, let us recall that, given $h\in\hat{\cal A}$, a weak solutions of the problem
$$\cases {-\left ( a+b\int_{\Omega}|\nabla u(x)|^2dx\right )\Delta u =h(x,u) & in $\Omega$\cr
& \cr u=0 & on $\partial\Omega$\cr}$$
is any $u\in H^1_0(\Omega)$ such that
$$\left ( a+b\int_{\Omega}|\nabla u(x)|^2dx\right )\int_{\Omega}\nabla u(x)\nabla w(x)dx=
\int_{\Omega}h(x,u(x))w(x)$$
for all $w\in H^1_0(\Omega)$. Let $\Phi:H^1_0(\Omega)\to {\bf R}$ be the functional defined by
$$\Phi(u)={{a}\over {2}}\|u\|^2+{{b}\over {4}}\|u\|^4$$
for all $u\in H^1_0(\Omega)$.\par
\smallskip
Hence, the weak solutions of the problem are precisely the critical points in
$H^1_0(\Omega)$ of the functional $\Phi-I_h$.\par
\smallskip
As an application of Theorem 5.8, we now obtain
\medskip
THEOREM 9.3. - {\it Let $n\geq 4$, let $f\in {\cal A}$ and let $g\in\tilde{\cal A}$ be such
that the set
$$\left \{x\in \Omega : \sup_{\xi\in {\bf R}}|g(x,\xi)|>0\right\}$$
has a positive measure.\par
Then,  there exists $\lambda^*\geq 0$ such that, for each $\lambda>\lambda^*$ and each 
convex set $C\subseteq L^2(\Omega)$ whose closure in $L^2(\Omega)$
contains the set $\{G(\cdot,u(\cdot)) : u\in H^1_0(\Omega)\}$,
there exists $v^*\in C$ such
that the problem
$$\cases {-\left ( a+b\int_{\Omega}|\nabla u(x)|^2dx\right )\Delta u =f(x,u)+\lambda(G(x,u)-v^*(x))g(x,u) & in $\Omega$\cr
& \cr u=0 & on $\partial\Omega$\cr}$$
has at least three weak solutions, two of which are global minima in $H^1_0(\Omega)$ of the functional
$$u\to {{a}\over {2}}\int_{\Omega}|\nabla u(x)|^2dx+{{b}\over {4}}\left ( \int_{\Omega}|\nabla u(x)|^2dx\right ) ^2-\int_{\Omega}F(x,u(x))dx
-{{\lambda}\over {2}}\int_{\Omega}|G(x,u(x))-v^*(x)|^2dx\ .$$
Furthermore, if the functional
$$u\to {{a}\over {2}}\int_{\Omega}|\nabla u(x)|^2dx+{{b}\over {4}}\left ( \int_{\Omega}|\nabla u(x)|^2dx\right ) ^2-\int_{\Omega}F(x,u(x))dx$$
has at least two global minima in $H^1_0(\Omega)$ and the function $G(x,\cdot)$ is strictly monotone for all $x\in \Omega$,
then $\lambda^*=0$.}\par
PROOF. For each $\lambda\geq 0$,
$v\in L^2(\Omega)$, consider the function $h_{\lambda,v}:\Omega\times {\bf R}\to {\bf R}$
defined by 
$$h_{\lambda,v}(x,\xi)=f(x,\xi)+\lambda(G(x,\xi)-v(x))g(x,\xi)$$ 
for all $(x,\xi)\in\Omega\times {\bf R}$.
Clearly, the function $h_{\lambda,v}$ lies in
$\hat{\cal A}$ and
$$H_{\lambda,v}(x,\xi)=F(x,\xi)+{{\lambda}\over {2}}\left ( |G(x,\xi)-v(x)|^2-|v(x)|^2\right )\ .$$
 So, the weak solutions of the problem are precisely the critical points in $H^1_0(\Omega)$
of the functional $\Phi-I_{h_{\lambda,v}}$. Moreover, if $p\in \left] 0,{{n+2}\over {n-2}}\right [$ and $q\in
\left ] 0,{{2}\over {n-2}}\right [$ are such that $(9.6)$ and $(9.7)$ hold, 
 for some constant $c_{\lambda,v}$, we have
$$\int_{\Omega}|H_{\lambda,v}(x,u(x))|dx\leq c_{\lambda,v}\left ( \int_{\Omega}|u(x)|^{p+1}+
\int_{\Omega}|u(x)|^{2(q+1)}dx+1\right )$$
for all $u\in H^1_0(\Omega)$. Therefore, by the Sobolev embedding theorem, for a constant $\tilde c_{\lambda,v}$,
we have
$$\Phi(u)-I_{h_{\lambda,v}}(u)\geq {{b}\over {4}}\|u\|^4-\tilde c_{\lambda,v}(\|u\|^{p+1}+\|u\|^{2(q+1)}+1) \eqno{(9.8)}$$
for all $u\in H^1_0(\Omega)$. On the other hand, since $n\geq 4$, one has
$$\max\{p+1, 2(q+1)\}<{{2n}\over {n-2}}\leq 4\ .$$
Consequently, from $(9.8)$, we infer that
$$\lim_{\|u\|\to +\infty}(\Phi(u)-I_{h_{\lambda,v}}(u))=+\infty\ .\eqno{(9.9)}$$
Since the functional $\Phi-I_{h_{\lambda,v}}$ is sequentially weakly lower semicontinuous, by the Eberlein-Smulyan theorem and
by $(9.9)$, it follows that it is inf-weakly compact.\par
Now, we are going to apply Theorem 5.8 taking $X=H^1_0(\Omega)$ with
the weak topology and $\Lambda=Y=L^2(\Omega)$ with the strong topology, and $y_0=0$. The symbols $x$ and $\lambda$ appearing
in Theorem 5.8 are replaced by the symbols $u$ and $v$ respectively. 
 Also, we take
$$\varphi(w)={{1}\over {2}}\int_{\Omega}|w(x)|^2dx$$
for all $w\in L^2(\Omega)$. Clearly, $\varphi\in {\cal G}$. Furthermore, we take
$$\Psi(u,v)(x)=G(x,u(x))-v(x)$$
for all $u\in H^1_0(\Omega), v\in L^2(\Omega), x\in \Omega$. Clearly, $\Psi(u,v)\in L^2(\Omega)$,
 $\Psi(u,\cdot)$ is a homeomorphism, and we have
$$v_u(x)=G(x,u(x))\ .$$
We show that the map $u\to v_u$ is not constant in $H^1_0(\Omega)$. 
Set
$$A=\left \{x\in \Omega : \sup_{\xi\in {\bf R}}|g(x,\xi)|>0\right\}\ .$$
By assumption, meas$(A)>0$. Then, by the classical Scorza-Dragoni theorem ([2], Theorem 2.5.19),
there exists a compact set $K\subset A$, of positive measure, such that the restriction
of $G$ to $K\times {\bf R}$ is continuous. Fix a point $\tilde x\in K$ such that
the intersection of $K$ and any ball centered at $\tilde x$ has a positive measure.
Next, fix $\xi_1, \xi_2\in {\bf R}$
such that
$$G(\tilde x,\xi_1)<G(\tilde x,\xi_2)\ .$$
By continuity, there is a closed ball $B(\tilde x,r)$ such that
$$G(x,\xi_1)<G(x,\xi_2)$$
for all $x\in K\cap B(\tilde x,r)$. Finally, consider two functions $u_1, u_2\in H^1_0(\Omega)$
which are constant in $K\cap B(\tilde x,r)$. So, we have
$$G(x,u_1(x))<G(x,u_2(x))$$
for all $x\in K\cap B(\tilde x,r)$. Hence, as meas($K\cap B(\tilde x,r))>0$, we infer that
$v_{u_1}\neq v_{u_2}$, as claimed. As a consequence, $\Psi\in {\cal H}$.  Of course,
$\varphi(\Psi(u,\cdot))$  is continuous and convex for all $u\in X$. 
Finally, take
$$J=\Phi-I_f\ .$$
Clearly, $J\in {\cal M}$. So, for what seen above, all the assumptions of Theorem 5.8 are satisfied.
Consequently, if we take 
$$\lambda^*=\theta(\varphi, \Psi,J)\eqno{(9.10)}$$ and
fix $\lambda>\lambda^*$ and a convex set  $C\subseteq L^2(\Omega)$ whose closure in $L^2(\Omega)$
contains the set $\{G(\cdot,u(\cdot)) : u\in H^1_0(\Omega)\}$, there exists $v^*\in C$ such that
the functional $\Phi-I_{h_{\lambda,v^*}}$ has at least two global minima in $H^1_0(\Omega)$ which are, therefore,
weak solutions of the problem. To guarantee the existence of a third solution, 
  denote by $k$ the inverse of the restriction of the function $at+bt^3$ to $[0,+\infty[$. 
Let $T:X\to X$ be the operator defined by
$$T(w)=\cases {{{k(\|w\|)}\over {\|w\|}}w & if $w\neq 0$\cr & \cr
0 & if $w=0$\ ,\cr}$$
 Since $k$ is continuous
and $k(0)=0$, the operator $T$ is continuous in $X$. For each $u\in X\setminus \{0\}$,
we have
$$T(\Phi'(u))=T((a+b\|u\|^2)u)={{k((a+b\|u\|^2)\|u\|)}\over {(a+b\|u\|^2)\|u\|}}(a+b\|u\|^2)u={{\|u\|}\over {(a+b\|u\|^2)\|u\|}}(a+b\|u\|^2)u=u\ .$$
In other words, $T$ is a continuous inverse of $\Phi'$. Then, since $I_{h_{\lambda,v^*}}'$ is compact, the functional
$\Phi-I_{h_{\lambda,v^*}}$ satisfies the Palais-Smale condition ([29], Example 38.25) and hence the existence of a third critical
point of the same functional is assured by Corollary 1 of [9].\par
Finally, assume that the functional $\Phi-I_f$ has at least two global minima, say $\hat u_1, \hat u_2$. Then,
the set $D:=\{x\in \Omega : \hat u_1(x)\neq \hat u_2(x)\}$ has a positive measure. By assumption, we have
$$G(x,\hat u_1(x))\neq G(x,\hat u_2(x))$$
for all $x\in D$, and so $v_{\hat u_1}\neq v_{\hat u_2}$. Then, by definition, we have
$$0\leq\theta(\varphi,\Psi,J)\leq {{J(\hat u_1)-J(\hat u_2)}\over {\varphi(\Psi(\hat u_1,v_{\hat u_2}))}}=0$$
and so $\lambda^*=0$ in view of $(9.10)$.\hfill $\bigtriangleup$\par
\medskip
Notice the following simple corollary of Theorem 9.3:\par
\medskip
COROLLARY 9.2. - {\it Let $n\geq 4$ and let $p\in \left ] 0,{{n+2}\over {n-2}}\right [$.\par
Then, for each $\lambda>0$ large enough and for each convex set $C\subseteq L^2(\Omega)$ whose closure in $L^2(\Omega)$
contains $H^1_0(\Omega)$, there exists $v^*\in C$ such that the problem
$$\cases {-\left ( a+b\int_{\Omega}|\nabla u(x)|^2dx\right )\Delta u =|u|^{p-1}u+\lambda(u-v^*(x)) & in $\Omega$\cr
& \cr u=0 & on $\partial\Omega$\cr}$$
has at least three solutions, two of which are global minima in $H^1_0(\Omega)$ of the functional
$$u\to {{a}\over {2}}\int_{\Omega}|\nabla u(x)|^2dx+{{b}\over {4}}\left ( \int_{\Omega}|\nabla u(x)|^2dx\right ) ^2-{{1}\over
{p+1}}\int_{\Omega}|u(x)|^{p+1}dx
-{{\lambda}\over {2}}\int_{\Omega}|u(x)-v^*(x)|^2dx\ .$$}\par
\medskip
Among the other consequences of Theorem 9.3, we highlight the following\par
\medskip
THEOREM 9.4. - {\it Let $n\geq 4$, let $f\in {\cal A}$ and let $g\in\tilde{\cal A}$ be such that
the set
$$\left\{x\in \Omega : \sup_{\xi\in {\bf R}}F(x,\xi)>0\right\}$$
haa a positive measure.
 Moreover, assume that, for each $x\in \Omega$, $f(x,\cdot)$ is
odd, $g(x,\cdot)$ is even and $G(x,\cdot)$ is strictly monotone.\par
Then,  for each $\lambda>0$, there exists $\mu^*> 0$ such that, for each $\mu>\mu^*$ and for each convex set $C\subseteq L^2(\Omega)$ whose closure
in $L^2(\Omega)$ contains the set $\{G(\cdot,u(\cdot)) : u\in H^1_0(\Omega)\}$, 
there exists $v^*\in C$ such
that the problem
$$\cases {-\left ( a+b\int_{\Omega}|\nabla u(x)|^2dx\right )\Delta u =\mu f(x,u)-\lambda v^*(x)g(x,u) & in $\Omega$\cr
& \cr u=0 & on $\partial\Omega$\cr}$$
has at least three weak solutions, two of which are global minima in $H^1_0(\Omega)$ of the functional
$$u\to {{a}\over {2}}\int_{\Omega}|\nabla u(x)|^2dx+{{b}\over {4}}\left ( \int_{\Omega}|\nabla u(x)|^2dx\right ) ^2-\mu\int_{\Omega}F(x,u(x))dx
+\lambda\int_{\Omega}v^*(x)G(x,u(x))dx\ .$$}\par
\smallskip
PROOF. Set
$$D=\left \{x\in \Omega : \sup_{\xi\in {\bf R}}F(x,\xi)>0\right\}\ .$$
By assumption, meas$(D)>0$. Then, by  the Scorza-Dragoni theorem,
there exists a compact set $K\subset D$, of positive measure, such that the restriction
of $F$ to $K\times {\bf R}$ is continuous. Fix a point $\hat x\in K$ such that
the intersection of $K$ and any ball centered at $\hat x$ has a positive measure.
Choose $\hat \xi\in {\bf R}$ so that $F(\hat x,\hat\xi)>0$. By continuity, there
is $r>0$ such that
$$F(x,\hat\xi)>0$$
for all $x\in K\cap B(\hat x,r)$. Set
$$M=\sup_{(x,\xi)\in \Omega\times [-|\hat \xi|,|\hat \xi|]}|F(x,\xi)|\ .$$
Since $f\in {\cal A}$, we have $M<+\infty$. Next, choose an open set $\tilde\Omega$
such that
$$K\cap B(\hat x,r)\subset\tilde\Omega\subset\Omega$$
and
$$\hbox {\rm meas}(\tilde\Omega\setminus (K\cap B(\hat x,r))<{{\int_{K\cap B(\hat x,r)}F(x,\hat \xi)dx}\over {M}}\ .$$
Finally, choose a function $\tilde u\in H^1_0(\Omega)$ such that
$$\tilde u(x)=\hat\xi$$
for all $x\in K\cap B(x,r)$,
$$\tilde u(x)=0$$
for all $x\in \Omega\setminus\tilde\Omega$ and
$$|\tilde u(x)|\leq |\hat\xi|$$
for all $x\in\Omega$. Thus, we have
$$\int_{\Omega}F(x,\tilde u(x))dx=\int_{K\cap B(\hat x,r)}F(x,\hat \xi)dx+\int_{\tilde\Omega\setminus (K\cap B(\hat x,r)}F(x,\tilde u(x))dx$$
$$>\int_{K\cap B(\hat x,r)}F(x,\hat \xi)dx-M\hbox {\rm meas}(\tilde\Omega\setminus (K\cap B(\hat x,r))>0\ .$$
Now, fix any $\lambda>0$ and set
$$\mu^*={{\Phi(\tilde u)+{{\lambda}\over {2}}I_{Gg}(\tilde u)}\over {I_f(\tilde u)}}\ .$$
Fix $\mu>\mu^*$. Hence
$$\Phi(\tilde u)-\mu I_f(\tilde u)+{{\lambda}\over {2}}I_{Gg}(\tilde u)<0\ .$$
From this, we infer that the functional $\Phi-\mu I_f+{{\lambda}\over {2}}I_{Gg}$ possesses at least to global
minima since it is even. At this point, we can apply Theorem 9.3 to the functions $g$ and $\mu f-\lambda Gg$. Our current
conclusion follows from the one of Theorem 9.3
since we have $\lambda^*=0$ and hence we can take the same fixed $\lambda>0$.\hfill $\bigtriangleup$
\medskip
To state the last result, denote by $H^{-1}(\Omega)$ the dual of $H^1_0(\Omega)$.\par
\smallskip
If $n\geq 2$, we denote by ${\cal A}$ the class of all
Carath\'eodory functions $f:\Omega\times {\bf R}\to {\bf R}$ such that
$$\sup_{(x,\xi)\in \Omega\times {\bf R}}{{|f(x,\xi)|}\over
{1+|\xi|^q}}<+\infty\ ,$$
where  $0<q< {{n+2}\over {n-2}}$ if $n>2$ and $0<q<+\infty$ if
$n=2$. While, when $n=1$, we denote by ${\cal A}$  the class
of all Carath\'eodory functions $f:\Omega\times {\bf R}\to {\bf R}$ such
that, for each $r>0$, the function $x\to \sup_{|\xi|\leq r}|f(x,\xi)|$ belongs
to $L^{1}(\Omega)$.\par
\smallskip
Given $f\in {\cal A}$ and $\varphi\in H^{-1}(\Omega)$,  consider the following Dirichlet problem
$$\cases {-\Delta u= f(x,u) +\varphi
 & in
$\Omega$\cr & \cr u=0 & on
$\partial \Omega$\ .\cr}\eqno{(P_{f,\varphi})} $$
 Let us recall
that a weak solution
of $(P_{f,\varphi})$ is any $u\in H^1_0(\Omega)$ such that
 $$\int_{\Omega}\nabla u(x)\nabla v(x)dx
-\int_{\Omega}f(x,u(x))v(x)dx-\varphi(v)=0$$
for all $v\in H^1_0(\Omega)$.\par
\smallskip
Applying Theorem 5.10, we obtain
\medskip
THEOREM 9.5. - {\it Let $f\in {\cal A}$ be such that the set
$$\left\{x\in \Omega : \sup_{\xi\in {\bf R}}\int_0^{\xi}f(x,t)dt>0\right\}$$
has a postive measure and
$$\limsup_{|\xi|\to +\infty}{{\sup_{x\in\Omega}\int_0^{\xi}f(x,t)dt}\over {\xi^2}}\leq 0\ .$$
Then, for every $\lambda>0$ large enough and for every convex set $C\subset H^{-1}(\Omega)$ dense
in $H^{-1}(\Omega)$, there exists $\varphi\in C$ such that
the problem
$$\cases {-\Delta u= \lambda f(x,u) +\varphi
 & in
$\Omega$\cr & \cr u=0 & on
$\partial \Omega$\cr} $$
has at least three weak solutions, two of which are global minima in $H^1_0(\Omega)$ of the functional 
$$u\to \int_{\Omega}|\nabla u(x)|^2dx-\lambda\int_{\Omega}\left (\int_0^{u(x)}f(x,\xi)d\xi\right )dx-\varphi(u)\ .$$}\par
\bigskip
\bigskip
{\bf Acknowledgement.} The author has been supported by the Gruppo Nazionale per l'Analisi Matematica, la Probabilit\`a e 
le loro Applicazioni (GNAMPA) of the Istituto Nazionale di Alta Matematica (INdAM). \par
\bigskip
\bigskip
\bigskip
\bigskip
\centerline {\bf References}\par
\bigskip
\bigskip
\noindent
[1]\hskip 5pt G. CORDARO, {\it On a minimax problem of Ricceri}, J. Inequal. Appl., {\bf 6} (2001), 261-285.\par
\smallskip
\noindent
[2]\hskip 5pt Z. DENKOWSKI, S. MIG\'ORSKI and N. S. PAPAGEORGIOU,
{\it An Introduction to Nonlinear Analysis: Theory}, Kluwer Academic
Publishers, 2003.\par
\smallskip
\noindent
[3]\hskip 5pt P. DR\'ABEK, {\it On a maximum principle for weak solutions of some
quasi-linear elliptic equations}, Appl. Math. Letters, {\bf 22}
(2009), 1567-1570.\par
\smallskip
\noindent
[4]\hskip 5pt P. DR\'ABEK and J. HERN\'ANDEZ, {\it Existence and
uniqueness of positive solutions for some quasilinear elliptic problems},
Nonlinear Anal., {\bf 44} (2001), 189-204.\par
\smallskip
\noindent
[5]\hskip 5pt V. L. KLEE, {\it Convexity of Chebyshev sets}, Math. Ann.,
{\bf 142} (1960/1961), 292-304.\par
\smallskip
\noindent
[6]\hskip 5pt E. KLEIN and A. C. THOMPSON, {\it Theory of Correspondences}, John Wiley
$\&$ Sons, 1984.\par
\smallskip
\noindent
[7]\hskip 5pt S. J. N. MOSCONI,
{\it A differential characterisation of the minimax inequality}, J. Convex
Anal., {\bf 19} (2012), 185-199.\par
\smallskip
\noindent
[8]\hskip 5pt R. A. PLASTOCK, {\it Nonlinear Fredholm maps of index zero and
their singularities}, Proc. Amer. Math. Soc., {\bf 68} (1978), 317-322.\par
\smallskip
\noindent
[9]\hskip 5pt  P. PUCCI and J. SERRIN, {\it A mountain pass theorem},
J. Differential Equations, {\bf 60} (1985), 142-149.\par
\smallskip
\noindent
[10]\hskip 5pt B. RICCERI, {\it A general variational principle and
some of its applications}, J. Comput. Appl. Math., {\bf 113}
(2000), 401-410.\par
\smallskip
\noindent
[11]\hskip 5pt B. RICCERI, {\it The problem of minimizing locally a $C^2$ functional
around non-critical points is well-posed}, Proc. Amer. Math. Soc., {\bf 135}
(2007), 2187-2191.\par
\smallskip
\noindent
[12]\hskip 5pt B. RICCERI, {\it Well-posedness of constrained minimization
problems via saddle-points}, J. Global Optim., {\bf 40} (2008),
389-397.\par
\smallskip
\noindent
[13]\hskip 5pt B. RICCERI, {\it On a theory by Schechter and Tintarev},
Taiwanese J. Math., {\bf 12} (2008), 1303-1312.\par
\smallskip
\noindent
[14]\hskip 5pt B. RICCERI, {\it Multiplicity of global minima for
parametrized functions}, Rend. Lincei Mat. Appl., {\bf 21} (2010),
47-57.\par
\smallskip
\noindent
[15]\hskip 5pt B. RICCERI, {\it A strict minimax inequality criterion and some of its consequences}, Positivity, {\bf 16} (2012), 455-470.\par
\smallskip
\noindent
[16]\hskip 5pt B. RICCERI, {\it Fixed points of nonexpansive potential
operators in Hilbert spaces}, Fixed Point Theory Appl. {\bf 2012},  2012: 123.\par
\smallskip
\noindent
[17]\hskip 5pt B. RICCERI, {\it A range property related to non-expansive
operators}, Mathematika, {\bf 60} (2014), 232-236.\par
\smallskip
\noindent
[18]\hskip 5pt B. RICCERI, {\it Integral functionals on $L^p$-spaces:
infima over sub-level sets}, Numer. Funct. Anal. Optim., {\bf 35} (2014), 1197-1211.\par
\smallskip
\noindent
[19]\hskip 5pt B. RICCERI, {\it Singular points of non-monotone potential operators}, J. Nonlinear
Convex Anal., {\bf 16} (2015), 1123-1129.\par
\smallskip
\noindent
[20]\hskip 5pt B. RICCERI, {\it Energy functionals of Kirchhoff-type problems having multiple global minima}, Nonlinear Anal., {\bf 115} (2015),  
130-136.\par
\smallskip
\noindent
[21]\hskip 5pt B. RICCERI, {\it Miscellaneous applications of certain minimax theorems I}, Proc. Dynam. Systems Appl., {\bf 7} (2016),
198-202.\par
\smallskip
\noindent
[22]\hskip 5pt R. T. ROCKAFELLAR, {\it  Integral functionals, normal integrands and measurable selections}, Lecture Notes in Math., Vol. 543,  157-207,
Springer, Berlin, 1976.\par
\smallskip
\noindent
[23]\hskip 5pt R. S. SADYRKHANOV, {\it On infinite                
dimensional features of proper and closed mappings}, Proc. Amer.
Math. Soc., {\bf 98} (1986), 643--658.\par
\smallskip
\noindent
[24]\hskip 5pt J. SAINT RAYMOND, {\it On a minimax theorem}, Arch. Math.
(Basel), {\bf 74} (2000), 432-437.\par
\smallskip
\noindent
[25]\hskip 5pt M. SCHECHTER and K. TINTAREV, {\it Spherical maxima in
Hilbert space and semilinear elliptic eigenvalue problems}, Differential
Integral Equations, {\bf 3} (1990), 889-899.\par
\smallskip
\noindent
[26]\hskip 5pt J. L. V\'AZQUEZ, {\it A strong maximum principle for some
quasilinear elliptic equations}, Appl. Math. Optim., {\bf 12} (1984),
191-202. \par
\smallskip
\noindent
[27]\hskip 5pt C. Z\u{A}LINESCU, {\it Convex analysis in general vector spaces}, World Scientific, 2002.\par
\smallskip
\noindent
[28]\hskip 5pt E. ZEIDLER, {\it Nonlinear functional analysis and its
applications}, vol. I, Springer-Verlag, 1986.\par
\smallskip
\noindent
[29]\hskip 5pt E. ZEIDLER, {\it Nonlinear functional analysis and its
applications}, vol. III, Springer-Verlag, 1985.\par
\bigskip
\bigskip
\bigskip
\bigskip
Department of Mathematics and Computer Science\par
University of Catania\par
Viale A. Doria 6\par
95125 Catania, Italy\par
{\it e-mail address}: ricceri@dmi.unict.it

\bye